\documentclass{article}
\usepackage{graphicx} 
\usepackage[english]{babel}
\usepackage{amsfonts,amsthm,amsmath,amssymb,mathtools,nicefrac}
\usepackage{fullpage}
\usepackage{authblk}
\usepackage{enumitem}
\usepackage{todonotes}
\usepackage{hyperref}
\usepackage{bbm,stmaryrd}
\usepackage[most]{tcolorbox}
 \usepackage{xr}
\usepackage{imakeidx}
\newcommand{\By}[2]{\overset{\mbox{\tiny{#1}}}{#2}}
\newcommand{\ByRef}[2]{   \By{\eqref{#1}}{#2} }
\newcommand{\eqBy}[1]{    \By{#1}{=} }

\newcommand{\leBy}[1]{    \By{#1}{\le} }

\newcommand{\eqByRef}[1]{ \ByRef{#1}{=} }

\newcommand{\leByRef}[1]{ \ByRef{#1}{\le} }

\newcommand{\JUSTIFY}[1]{\fbox{\tiny{#1}}\quad}

\newtheorem{lem}{Lemma}[section]

\newtheorem{defi}[lem]{Definition}
\newtheorem{rem}[lem]{Remark}
\newtheorem{thm}[lem]{Theorem}

\newtheorem{prop}[lem]{Proposition}

\newcommand{\OriginalTwoSAT}{\mathsf{TwoSAT}}
\newcommand{\DensestTwoSAT}{\mathsf{DensestTwoSAT}}
\newcommand{\OriginalKSAT}{\mathsf{Random}$k$\mathsf{SAT}}
\newcommand{\TwoSAT}{\mathsf{TwoSAT}^\dagger}
\newcommand{\Lit}{\mathrm{Lit}}
\newcommand{\PSpec}{\mathrm{PSpec}}
\newcommand{\Spec}{\mathrm{Spec}}

\newcommand{\OurSpace}{\mathbf{K}}
\newcommand{\OurMeasure}{\boldsymbol{\kappa}}
\newcommand{\G}{\mathbb{G}}
\newcommand{\N}{\mathbb{N}}

\newcommand{\R}{\mathbb{R}}

\newcommand{\Probability}{\mathbf{P}} 
\newcommand{\Expectation}{\mathbf{E}} 
\newcommand{\diff}{\mathsf{d}}
\newcommand{\eps}{\varepsilon}
\newcommand{\signum}{\mathfrak{S}}
\newcommand{\ImaginaryUnit}{\mathbf{i}}
\newcommand{\ImplDig}{\mathsf{D}} 

\newcommand{\BoundOper}{\mathcal{B}}
\newcommand{\CompactOper}{\mathcal{K}}

\newcommand{\Krein}{Kre\u{\i}n}

\DeclareMathOperator{\esssup}{ess\, sup}
\DeclareMathOperator{\essinf}{ess\, inf}

\newcommand{\impl}{\mathrm{impl}}
\newcommand{\Var}{\mathrm{Vars}}
\newcommand{\clauses}{\mathrm{clauses}}

\title{Inhomogeneous random 2-SAT}
\author{Jan Hladký, Petr Savický}
\date{}
\affil{Institute of Computer Science of the Czech Academy of Sciences\thanks{Research supported by Czech Science Foundation Project 21-21762X. With institutional support RVO:67985807.}}

\makeindex

\begin{document}
\maketitle
\begin{abstract}
We introduce an inhomogeneous variant of random 2-SAT. Each variable $v_1,\ldots,v_n$ is assigned a type from a state space $\Lambda$, independently at random. Clause inclusion is governed by a symmetric measurable kernel $W$ on $(\Lambda \times \{+,-\})^2$, in analogy with the inhomogeneous random graph model of Bollobás, Janson, and Riordan: given literals $\ell_i\in\{v_i,\neg v_i\}$ and $\ell_j\in\{v_j,\neg v_j\}$, the clause $\{\ell_i,\ell_j\}$ appears with probability $W(\mathrm{type}(\ell_i),\mathrm{type}(\ell_j))/(2n)$. In particular, for a variable $v_i$ of type $x\in\Lambda$, the slices $W((+,x),\cdot)$ and $W((-,x),\cdot)$ describe how $v_i$ and $\neg v_i$ interact with other literals.

We identify a parameter $\rho^*(W)$, defined as the spectral radius of an integral operator derived from $W$, and show that $\rho^*(W)<1$ and $\rho^*(W)>1$ correspond to asymptotically almost surely satisfiable and unsatisfiable instances, respectively. 
The satisfiability threshold for homogeneous random 2-SAT is well-established, occurring when the ratio of clauses to variables is $1$. This corresponds to a weight function of $W \equiv 1$ and a clause density of $1/(2n)$.
Our result extends this classical result to a broad class of models controlled by types of variables.
\end{abstract}

\newpage
{
\renewcommand{\baselinestretch}{1}\normalsize
\tableofcontents
}
\newpage

\section{Introduction}\label{sec:Intro}

For $k\in\{2,3,\ldots\}$, \emph{random $k$-SAT} is one of the most studied problems in the intersection of combinatorics, probability theory, and theoretical computer science. It is a natural probabilistic model in the area of satisfiability, where the goal is to determine whether a given Boolean formula can be satisfied---i.e., whether there exists an assignment of truth values to variables that makes the entire formula evaluate to true.

In the $k$-SAT problem, the formula is expressed in conjunctive normal form and consists of variables $v_1,\ldots,v_n$ and $m$ clauses, where each clause is a disjunction of $k$ literals. A literal is either a variable or its negation. A random $k$-SAT instance is generated by selecting each clause uniformly at random from the set of all possible clauses involving $k$ distinct literals.

One of the key aspects of interest in random $k$-SAT is the behavior of the problem as the ratio of the number of clauses $m$ to the number of variables $n$, denoted by $\alpha = \frac{m}{n}$, varies. As $\alpha$ increases, it is believed that the probability of the formula being satisfiable undergoes a sharp transition, commonly referred to as the \emph{satisfiability threshold $\alpha_k$}. Conjecturally, for values of $\alpha$ below $\alpha_k$, a random $k$-SAT for $\alpha$ is satisfiable asymptotically almost surely (as $n\to\infty$), while for $\alpha$ above it, it becomes unsatisfiable asymptotically almost surely. As an alternative to parametrizing the model by $n$ and $m$, we can parametrize the model by $n$ and $\alpha$. In this latter model, we form a formula by including each possible clause (there are $2^k\binom{n}{k}$ such clauses) with probability $\alpha n 2^{-k}\binom{n}{k}^{-1}$. So, the difference between the `$(n,m)$-model' and the `$(n,\alpha)$-model' of random $k$-SAT is the same as the difference between the uniform model and the binomial model of Erdős--Rényi random graphs, and is insignificant for our purposes.

It is well-known that the $k$-SAT is in the computational complexity class $\mathsf{P}$ for $k=2$ whereas it is $\mathsf{NP}$-complete for every $k>2$. In fact, it is known that each 2-SAT formula can be represented by an `implication digraph' and that satisfiability of the formula then corresponds to a simple and computationally tractable problem of the existence of certain `contradictory cycles'. There is a similar jump in the arduousness of the analysis of random $k$-SAT. In particular, in 1992, Goerdt~\cite{Goerdt92} and independently Chvátal and Reed~\cite{ChvatalR92} proved the above `satisfiability conjecture', for $k=2$ and further determined that $\alpha_2=1$. Using the connection above, the bulk of their proof is in fact about random digraphs. On the other hand, for every $k\ge 3$, even the existence of the satisfiability threshold $\alpha_k$ was unknown (albeit Friedgut~\cite{MR1678031} got close). This changed only in~2014 with a preprint version of a \emph{tour de force} paper~\cite{MR4429261}, which used methods of statistical physics to determine~$\alpha_k$ for all $k$ sufficiently large. Among a huge body of work which looks at various aspects of random $k$-SAT, we chose~\cite{MR3161470,MR3436404,MR2083472,MR4258129,MR1824274} as representative references. We also refer to a slightly outdated survey~\cite{SurveyAchlioptas09}.

While random $k$-SAT is an elegant mathematical model, many real-world scenarios that can be effectively represented by random models often involve constraints that are not uniformly distributed. A prominent example in this direction is random 2-SAT with prescribed literal degrees studied in~\cite{MR2334579}. For each $n$, this model is parametrized by integers $d_{1},\bar{d_{1}},d_{2},\bar{d_{2}},\ldots,d_{n},\bar{d_{n}}\in\N_0$ with $\sum_i (d_{i}+\bar{d_{i}})$ even. Further, it is assumed that $\max_i\{d_{i},\bar{d_{i}}\}\le n^{1/11}$. A random formula is then taken uniformly from the set of all simple 2-SAT formulae\footnote{a 2-SAT formula is \emph{simple} if it does not contain a repetition of any clauses, and also the literals of each clause arise from different variables} with exactly $d_{i}$ occurrences of literal $v_i$ and exactly $\bar{d_{i}}$ occurrences of literal $\neg v_i$. The main result of~\cite{MR2334579} asserts that depending on whether $2\sum_id_{i}\cdot \bar{d_{i}}<(1-\eps)\sum_i(d_{i}+\bar{d_{i}})$ or $2\sum_id_{i}\cdot \bar{d_{i}}>(1+\eps)\sum_i(d_{i}+\bar{d_{i}})$, a random 2-SAT formula with prescribed literal degrees is asymptotically almost surely satisfiable or asymptotically almost surely unsatisfiable, respectively.

\subsection{Our model $\OriginalTwoSAT(n,W)$}\label{ssec:OriginalSat}
Whenever we refer to a subset of a measure space, we implicitly assume it is measurable. Suppose that $B$ is a probability space with measure $\beta$. Suppose that $p\in[1,\infty]$. Nonnegative functions in $L^p(B^2)$ are called
\emph{$L^p$-digraphons}. $L^\infty$-digraphons are simply called \emph{digraphons}. $L^p$-digraphons which are symmetric with respect to the swap of the coordinates are called \emph{$L^p$-graphons}. Again, $L^\infty$-graphons are simply called \emph{graphons}. Note that compared to other literature, the values in graphons or digraphons in this paper are not necessarily bounded from above by~1. The combinatorial interpretation of the traditional bound of~1 is that each pair of vertices is connected by at most~1 edge. In this paper, however, we use graphons and digraphons mostly as sources of sparsified random graphs, that is, for a digraphon $W$, the edge inclusion probabilities are encoded in $\frac{W}{n}$, where $n\in \N$ is large.

Throughout the paper, \index{$(\Lambda,\lambda)$}$\Lambda$ is an arbitrary Polish space and $\lambda$ is a Borel probability measure on it. Let
\index{$(\signum,\mu^{+-})$}$\signum:=\{+,-\}$ be the two-element probability space equipped with the
uniform measure $\mu^{+-}$. The space $\signum$ represents the positive and negative signs we equip the logical variables in our formula with. Let us consider the probability space \index{$(\OurSpace,\OurMeasure)$}$\OurSpace=\Lambda\times\signum$ with the corresponding product measure $\OurMeasure=\lambda\times \mu^{+-}$.
Define the negation map \index{$\neg$}$\neg:\OurSpace \to \OurSpace$
by $\neg (x,+) = (x,-)$ and $\neg (x,-) = (x,+)$ for every $x \in \Lambda$.
Let $W$ be an  $L^1$-graphon on $\OurSpace$. For $n\in \N$, we define random 2-SAT
formula \index{$\OriginalTwoSAT(n,W)$}$\phi\sim\OriginalTwoSAT(n,W)$ on variables $\{v_1,\ldots,v_n\}$
and let us denote the set of literals \index{$\Lit_n$}$\Lit_n = \{v_1,\ldots,v_n,\neg v_1, \ldots, \neg v_n\}$. For this purpose, we sample elements
$x_1,\ldots,x_n\in \Lambda$ independently with distribution $\lambda$. For each
pair $\{i,j\}\in\binom{n}{2}$ we insert clauses
$\{v_i, v_j\}$, $\{v_i, \neg v_j\}$, $\{\neg v_i, v_j\}$, $\{\neg v_i, \neg v_j\}$
into $\phi$ independently at random with probabilities specified as follows.
For $\mathfrak{q},\mathfrak{s}\in\signum$, the clause
$\{\mathfrak{q} v_i , \mathfrak{s} v_j\}$ is inserted with probability
\begin{equation}\label{eq:def2SATjb}
\min\left\{1,\frac{W\left((x_i,\mathfrak{q}),(x_j,\mathfrak{s})\right)}{2n}\right\}    
\end{equation}
where we identify $+v$ with $v$ and $-v$ with $\neg v$ for every logical variable $v$.
Also, note that we view the clauses as unordered pairs, so that $\{\mathfrak{q} v_i , \mathfrak{s} v_j\}=\{\mathfrak{s} v_j , \mathfrak{q} v_i\}$. Therefore, the fact that $W$ is symmetric is used in that~\eqref{eq:def2SATjb} does not depend on the order of the literals.

\subsubsection{Stochastic block model}\label{ssec:stochblockmodel}
While the number of clauses in $\OriginalTwoSAT(n,W)$ is random, it is easy to see that it is concentrated (as $n\to\infty$) at $(1\pm o(1)) n \|W\|_1$. Note that the original  $(n,\alpha)$-model of random 2-SAT corresponds to $\OriginalTwoSAT(n,\boldsymbol{\alpha})$, where $\boldsymbol{\alpha}$ is the constant-$\alpha$ function. A class of models between the original random 2-SAT and the full generality of $\OriginalTwoSAT(n,W)$ could be called `stochastic block model'.\footnote{This term is borrowed from a related model of random graphs that first appeared in relation to a problem in sociology~\cite{MR718088} and has been widely used since.} This model is parametrized by the \emph{number of types} $t\in \N$, the \emph{proportions} $\gamma_1,\ldots,\gamma_t>0$, $\sum \gamma_i=1$, and \emph{connection parameters} $C^{\mathfrak{q},\mathfrak{s}}_{i,j}\ge 0$, $\mathfrak{q},\mathfrak{s}\in\signum$, $i,j\in [t]$ subject to symmetry $C^{\mathfrak{q},\mathfrak{s}}_{i,j}=C^{\mathfrak{s},\mathfrak{q}}_{j,i}$. For a given $n$, we take numbers $n_1,\ldots,n_t$, where $n=n_1+\ldots+n_t$, $\frac{n_i}{n}\approx \gamma_i$ for each $i\in[t]$ (subject to some approximation condition) and Boolean variables $\{v_{i,k}\}_{i\in[t],k\in[n_i]}$. Each clause
$\{\mathfrak{q} v_{i,k} , \mathfrak{s} v_{j,\ell}\}$ on distinct variables is included with probability $C^{\mathfrak{q},\mathfrak{s}}_{i,j}/(2n)$. This, for many practical purposes, corresponds to $\OriginalTwoSAT(n,W)$ if $W$ is a graphon defined on $\Lambda\times\signum$, where $\Lambda=\Lambda_1\sqcup \ldots\sqcup \Lambda_t$, $\lambda(\Lambda_i)=\gamma_i$ and $W_{\restriction (\Lambda_i\times \{\mathfrak{q}\})\times (\Lambda_j\times \{\mathfrak{s}\})}\equiv C^{\mathfrak{q},\mathfrak{s}}_{i,j}$. To see this, consider the stage of sampling elements $x_1,\ldots,x_n\in \Lambda$ in the procedure of generating $\OriginalTwoSAT(n,W)$. By the Law of Large Numbers, the number of indices $\ell\in [n]$ for which $x_\ell\in \Lambda_i$ (for a given $i\in[t]$) satisfies with high probability that $\frac{n_i}{n}\approx \gamma_i$. Given this event, the individual clauses are inserted with the same probability as in the stochastic block model.

\subsection{The statement of the result}
Our main result, Theorem~\ref{thm:main}, demonstrates that the threshold
phenomenon known for the homogeneous random 2-SAT also extends to $\OriginalTwoSAT(n,W)$. While the model $\OriginalTwoSAT(n,W)$ is sensible for any nonnegative symmetric measurable function $W$ on $\OurSpace^2$, in our main theorem we will impose mild additional integrability and operator-boundedness conditions. The critical parameter is characterized by the spectral properties of a specific operator on the Banach space $L^1(\OurSpace)$ derived from the $L^1$-graphon $W$. We introduce several concepts to state the result. In Definition~\ref{def:implicationdigraphon} we introduce the implication $L^p$-digraphon $\overrightarrow{W}$. In Definition~\ref{def:restriction} we introduce a restriction of an $L^p$-digraphon to a set. In Definition~\ref{def:component} we introduce strong components of an $L^1$-digraphon and in Theorem~\ref{thm:decompositionIntoComponents} we state the existence of decomposition of an $L^1$-digraphon into its strong components, paralleling classical results on digraphs. In Definition~\ref{def:contradictoryset} we define contradictory sets. Finally, in Definition~\ref{def:spectrum} we recall the notions of eigenvalues and spectral radius of $L^p$-digraphons. So, while the central theme of the paper is random 2-SAT, some of the tools we develop have broader applications and are separated for clarity into \cite{HladkySavicky:Digraphons}: the decomposition of digraphons into strong components and the study of their spectral properties.

\begin{figure}\centering
	\includegraphics[scale=0.7]{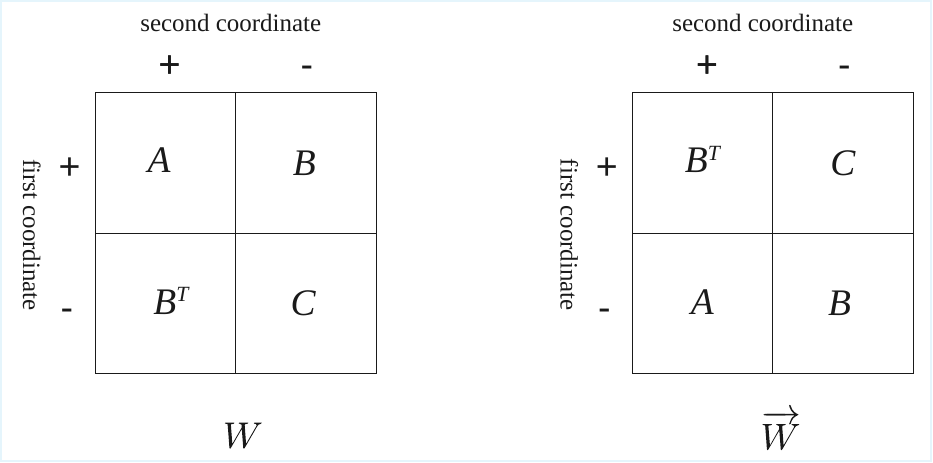}
	\caption{A visualization of the transformation of a graphon $W$ into its implication digraphon $\protect\overrightarrow{W}$. The graphon $W$ consists of four parts: a symmetric part $A\in L^1\left((\Lambda\times\{+\})\times (\Lambda\times\{+\})\right)$, a symmetric part $C\in L^1\left((\Lambda\times\{-\})\times (\Lambda\times\{-\})\right)$ and a pair of mutually transposed parts $B$ and $B^T$, where $B\in L^1\left((\Lambda\times\{+\})\times (\Lambda\times\{-\})\right)$, $B^T\in L^1\left((\Lambda\times\{-\})\times (\Lambda\times\{+\})\right)$. (The orientation of the plane here and elsewhere follows the matrix convention, that is, the main diagonal is in the $\searrow$ direction.)}
	\label{fig:SATKernelTransposition}
\end{figure}

In Section~\ref{sec:Intro} we mentioned that a 2-SAT formula can be turned into an `implication digraph'. We recall this well-known transformation in Section~\ref{sec:2SATUsingImplicationDigraph}. This transformation inspires a transformation of a graphon $W$ (which parametrizes $\OriginalTwoSAT(n,W)$) into an `implication digraphon' $\overrightarrow{W}$. See Figure~\ref{fig:SATKernelTransposition} for an illustration.
\begin{defi}[implication digraphon]\label{def:implicationdigraphon}
Suppose that $p\ge 1$ is given.
Given an $L^p$-graphon $W$ on $\OurSpace=\Lambda\times\signum$, its
\emph{implication $L^p$-digraphon} \index{$\overrightarrow{W}$}$\overrightarrow{W}$ is an $L^p$-digraphon on $\OurSpace$ defined for $x,y\in \OurSpace$ by $\overrightarrow{W}\big(x,y\big) := W\big(\neg x,y\big)$. 
\end{defi}
Note that the property of symmetry is in general lost in the implication digraphon. By using Definition~\ref{def:implicationdigraphon}, symmetry of $W$, and then Definition~\ref{def:implicationdigraphon} again, we have 
\begin{equation}\label{eq:antisymmetric}
\overrightarrow{W}\big(x,y\big) := W\big(\neg x,y\big)=W\big(y,\neg x\big)=\overrightarrow{W}\big(\neg y,\neg x\big)
\end{equation}
In the case of graphs, a corresponding property is called skew-symmetry.

Next, we define restriction of a digraphon to a subset of its ground set.
\begin{defi}[restriction]\label{def:restriction}\index{$U\llbracket A\rrbracket$}
Suppose that $p\ge 1$ is given.
For an $L^p$-digraphon $U$ on a probability space $(\Omega,\mu)$, and for $A\subseteq \Omega$, let $U\llbracket A\rrbracket$ be an $L^p$-digraphon on $\Omega$ defined by 
\[
U\llbracket A\rrbracket(x,y)=\begin{cases}U(x,y) &\mbox{if $x,y\in A$}\\
0&\mbox{otherwise.}
\end{cases}
\]
\end{defi}

In Definition~\ref{def:component} below, we introduce connectivity notions for digraphons taken from~\cite{HladkySavicky:Digraphons}. To motivate them by finite graphs, recall that a nonempty set of vertices $C$ in a digraph is \emph{strongly connected} if for every partition $A\sqcup B=C$ into two nonempty sets, there is at least one directed edge going from $A$ to $B$. Maximal strongly connected sets are called \emph{strong components}. These notions have almost straightforward counterparts for digraphons except they are uniquely defined only modulo nullsets.
The last notion in Definition~\ref{def:component} is of `fragmented sets'. This concept does not make any sense in finite digraphs. That is, one way to decompose a digraph into strong components is to start with an initial one-cell partition of the entire vertex set. If, at any moment, there is a cell violating the above condition on strong connectedness, we subdivide that cell accordingly. In a finite graph, this process must eventually terminate (since single vertices cannot be subdivided), resulting in the unique decomposition into strong components. In a digraphon, the regions where this process continues to split sets of positive measure into ever smaller ones constitute its fragmented sets.

\begin{defi}[strongly connected set, strong component, fragmented set]\label{def:component}
Suppose that $\Gamma$ is an $L^1$-digraphon on a probability space $(\Omega,\mu)$, and $X\subseteq \Omega$ is a set of positive measure.
\begin{enumerate}[label=(\roman*)]
    \item\label{en:StrongConn} We say that \emph{$X$ is strongly connected in $\Gamma$} if for every partition $A\sqcup B=X$ with $\mu(A),\mu(B)>0$ we have $\int_{A\times B}\Gamma>0$.
    \item\label{en:StrongComp} We say that \emph{$X$ is a strong component in $\Gamma$} if $X$ is strongly connected and for every $Y\subset \Omega$ with $\mu(X\cap Y)>0$ and $\mu(Y\setminus X)>0$ we have that $Y$ is not strongly connected.
    \item\label{en:Fragmented} We say that \emph{$X$ is fragmented in $\Gamma$} if every subset $Y\subset X$ of positive measure is not strongly connected.
\end{enumerate}
\end{defi}
Recall that each digraph can be decomposed in a unique way into maximal strong components. Theorem~\ref{thm:decompositionIntoComponents} below, which is one of the main results of \cite{HladkySavicky:Digraphons}, is a digraphon counterpart.\footnote{We state Theorem~\ref{thm:decompositionIntoComponents} in the generality of $L^1$-digraphons whereas the original statement in~\cite{HladkySavicky:Digraphons} is only about digraphons. Given a general $L^1$-digraphon $\Gamma$, we can, however, consider its \emph{indicator digraphon} $\Gamma'(x,y)=\mathbbm{1}_{\Gamma(x,y)>0}$. It is easy to see that the notions of strong components and fragmented sets for $\Gamma$ and for $\Gamma'$ are equivalent, and thus we have reduced the general case of $L^1$-digraphons to that of digraphons.}

\begin{thm}[Theorem~\ref{DIGRAPHONS-thm:decompositionIntoComponents} in~\cite{HladkySavicky:Digraphons}]\label{thm:decompositionIntoComponents}
Suppose that $\Gamma$ is an $L^1$-digraphon on $\Omega$. Then there exists a finite or a countable set $I$ not containing~$0$ and a decomposition $\Omega=X_0\sqcup \bigsqcup_{i\in I} X_i$ so that $X_0$ is either an empty set or is fragmented in $\Gamma$ and each $X_i$ is a strong component. 

Further, this decomposition is unique in the sense that if partitions $\{X_i\}_{i\in I\cup\{0\}}$ and $\{X'_i\}_{i\in I'\cup\{0\}}$ are two decompositions of $\Gamma$ into strong components as above, then there exists a bijection $\pi:I\to I'$ such that $X'_{\pi(i)}$ equals $X_i$ modulo a nullset for each $i\in I$, and $X_0$ equals $X'_0$ modulo a nullset.
\end{thm}

In Section~\ref{sec:Intro} we mentioned that the analysis of random 2-SAT goes via translating a 2-SAT formula into its implication digraph, and its satisfiability boils down to the existence of contradictory cycles. We will introduce these concepts in Section~\ref{sec:2SATUsingImplicationDigraph}. In our main theorem, Theorem~\ref{thm:main}, we determine which models $\OriginalTwoSAT(n,W)$ yield almost surely satisfiable or almost surely unsatisfiable formulae. To this end, we decompose the implication digraphon $\overrightarrow{W}$ into strong connected components and disregard those components that cannot generate contradictory cycles. To this end, we use the concept of contradictory sets below.
\begin{defi}[contradictory set]\label{def:contradictoryset}
A set $X\subseteq \OurSpace$ is \emph{contradictory} if $\lambda(\{x\in\Lambda:(x,+),(x,-)\in X\})>0$. 
\end{defi}

Spectral properties of digraphons are essential for formulating and proving our main result. While the spectral theory for graphons is well developed (see, e.g., Section~11.6 in~\cite{MR3012035}), its extension to digraphons poses nontrivial challenges. The main reason why spectral theory does not transfer directly is that graphons correspond to self-adjoint operators and thus admit a standard spectral decomposition, whereas digraphons generally do not. At the outset of this project, the spectral theory for digraphons was largely undeveloped. Addressing this gap --- and, more broadly, developing basic tools related to spectral analysis and connectivity  --- led us to write a separate, self-contained paper~\cite{HladkySavicky:Digraphons}. In the meantime, Grebík, Král', Liu, Pikhurko, and Slipantschuk~\cite{grebik2025convergencespectradigraphlimits} posted a preprint which studies spectral properties of digraphons in the context of generating random oriented graphs, including an expression for the homomorphism density of oriented cycles in terms of the point spectrum. For now, we introduce only the essential definitions necessary to state Theorem~\ref{thm:main}, and recall the spectral tools from~\cite{HladkySavicky:Digraphons} in Section~\ref{ssec:SpectralTheory}.

We fix $p\in [1,\infty)$. We work in the complex Banach space $L^p(\Omega)$, where $\Omega$ is a measure space (with an implicit sigma-algebra) equipped with a probability measure $\mu$.  For every $L^1$-digraphon $\Gamma$ on $\Omega$, we may consider integral kernel operator $T_\Gamma:L^p(\Omega)\to L^p(\Omega)$, where for $f\in L^p(\Omega)$ the function $g:=T_\Gamma(f)$ is defined by
\[
g(x):=\int_{y\in\Omega}f(y)\Gamma(y,x) \diff\mu(y)\;\mbox{for every $x\in \Omega$}.
\]
Note that given $p\in [1,\infty)$, it is not automatic that $g\in L^p(\Omega)$, or even that the integral to define $g(x)$ is finite. Only when this is the case, we can view $T_\Gamma$ as an integral kernel operator on $L^p(\Omega)$, and write \index{$\BoundOper(L^p(\Omega))$}$\Gamma\in\BoundOper(L^p(\Omega))$.
Henceforth, any reference to a spectral concept of an $L^1$-digraphon $\Gamma$ will mean the corresponding concept for its associated integral kernel operator $T_\Gamma$ (viewed on a Banach space $L^p(\Omega)$ for $p\in[1,\infty)$ which will be specified). We also write $\Gamma$ instead of $T_\Gamma$.

\begin{defi}[eigenvalues, eigenfunctions, spectrum, point spectrum, spectral radius]\label{def:spectrum}
Suppose that $p\in [1,\infty)$. Suppose that $\Gamma$ is an $L^1$-digraphon on $\Omega$ with the property that $\Gamma\in\BoundOper(L^p(\Omega))$. A complex number $\tau$ and a nonzero function $f\in L^p(\Omega)$ are called \emph{eigenvalue} and \emph{eigenfunction} of $\Gamma$ with respect to the Banach space $L^p(\Omega)$, respectively, if $\Gamma(f)=\tau f$. The collection of all eigenvalues and $0$ if $\Gamma$ is not invertible is called the \emph{point spectrum} of $\Gamma$ with respect to the Banach space $L^p(\Omega)$, and denoted \index{$\PSpec_p(\Gamma)$}$\PSpec_p(\Gamma)$. 
The collection of all complex numbers $\tau$ for which $\Gamma-\tau\cdot \mathbbm{1}$ is not invertible (as a bounded operator), is called the \emph{spectrum} of $\Gamma$, and denoted by $\Spec_p(\Gamma)$. 
The \emph{spectral radius} of $\Gamma$, denoted $\rho_p(\Gamma)$\index{$\rho_p(\cdot)$}, is the supremum of moduli taken over the elements of the spectrum. For the most important choice $p=2$, we write \index{$\rho(\cdot)$}$\rho(\cdot):=\rho_2(\cdot)$.
\end{defi}
The definitions and the properties implicitly implied in Definition~\ref{def:spectrum} are standard, see P7.3-5 in~\cite{MR0467220}. It is well-known that if in addition $\Gamma\in\CompactOper(L^p(\Omega))$, we have $\Spec_p(\Gamma)=\PSpec_p(\Gamma)$. Therefore, the spectral radius is the maximum modulus of an eigenvalue of $\Gamma$.

\begin{rem}\label{rem:quiteoften}
Quite often, $\rho_p(\Gamma)$ does not depend on $p$, for a wide range of values of $p$. As a prominent example, it follows from Lemma~\ref{lem:LpImpliesBoundedCompact} and Lemma~\ref{lem:GelfandInvariant}\ref{en:lGI2} that if $\Gamma\in L^p(\Omega^2)$ for some $p\in [2,\infty)$, then we have $\Gamma\in\bigcap_{q\in[2,p]}\CompactOper(L^q(\Omega))$ and for all $q\in [2,p]$ we have $\rho_p(\Gamma)=\rho_q(\Gamma)$. As another example which follows from Lemma~\ref{lem:eigenvectorsbounded}, when $\Gamma$ is a digraphon, then we have $\rho(\Gamma)=\rho_p(\Gamma)$ for all $p\in [1,\infty)$.
\end{rem}

Our main result, Theorem~\ref{thm:main} below, asserts that the asymptotic almost sure (un)satisfiability is determined by a certain spectral parameter $\rho^*_p(W)$. We define this parameter here.
\begin{defi}\label{defi:rho*}
Let $W$ be an $L^1$-graphon on $\OurSpace=\Lambda\times\signum$. Let $\OurSpace=\Omega_0\sqcup \bigsqcup_{i\in I}\Omega_i$ be a decomposition of $\overrightarrow{W}$ into strong components. Let $I^*\subseteq I$ be the indices of the contradictory components $\Omega_i$. For any $p\in[1,\infty)$ for which $\overrightarrow{W}\left\llbracket \bigcup_{i\in I^*}\Omega_i\right\rrbracket\in\CompactOper(L^p(\Omega))$, define \index{$\rho^*_p(\cdot)$}$\rho^*_p(W):=\sup_{i\in I^*}\rho_p\left(\overrightarrow{W}\llbracket \Omega_i\rrbracket\right)$. For the most important choice $p=2$, we write \index{$\rho^*(\cdot)$}$\rho^*(\cdot):=\rho^*_2(\cdot)$.
\end{defi}
We can now state the main result, Theorem~\ref{thm:main}. The operator-theoretic assumption on the $L^1$-graphon $W$ involved in Part~\ref{en:MainUnsatisfiable} may look a bit technical. As we noted in Remark~\ref{rem:quiteoften}, it is satisfied when $W\in L^2(\OurSpace^2)$. In particular, when $W\in L^3(\OurSpace^2)$, both the satisfiability and the unsatisfiability part are applicable.
\begin{thm}\label{thm:main}
    Let $W$ be an $L^1$-graphon on $\OurSpace=\Lambda\times\signum$.
    \begin{enumerate}[label=(\roman*)]
        \item\label{en:MainSatisfiable} If $W\in L^3(\OurSpace^2)$ and $\rho^*_3(W)<1$, then as $n\to\infty$, $\OriginalTwoSAT(n,W)$ is asymptotically almost surely satisfiable.
        \item\label{en:MainUnsatisfiable} Suppose that there exists $p\in [1,\infty)$ for which $\overrightarrow{W}\in\CompactOper(L^p(\Omega))$, and $\rho^*_p(W)>1$. Then, as $n\to\infty$, $\OriginalTwoSAT(n,W)$ is asymptotically almost surely unsatisfiable.
    \end{enumerate}
    \end{thm}
We prove Theorem~\ref{thm:main}\ref{en:MainSatisfiable} in Section~\ref{sec:ProofSatisfiable} and Theorem~\ref{thm:main}\ref{en:MainUnsatisfiable} in Section~\ref{sec:Unsatisfiable}. In Section~\ref{ssec:Heuristics} we give a basic idea behind the proof.

\begin{rem}\label{rem:OnIntegrabilityConditions}
In many applications, including the stochastic block model from Section~\ref{ssec:stochblockmodel}, the $L^1$-graphon $W$ is bounded in $L^\infty$. In that case, $\rho^*_p(W)$ in Definition~\ref{defi:rho*} is defined for every $p\in[1,\infty)$ and does not depend on the choice of $p$. Other prominent examples where $\rho^*_p(W)$ stays constant for many choices of $p$ were given in Remark~\ref{rem:quiteoften}.

Most prominent applications of unbounded graphons involve random scale-free formulea. These are discussed in Section~\ref{ssec:scalefree}. 

We believe that the integrability assumption $W\in L^3(\OurSpace^2)$ in Theorem~\ref{thm:main}\ref{en:MainSatisfiable} is an artifact of our proof and that it can be relaxed (at least) to $W\in L^2(\OurSpace^2)$.
\end{rem}


\begin{rem}
Theorem~\ref{thm:main} does not cover the case $\rho^*(W)=1$. In such a case, since the parameter $\rho^*(W)$ is easily seen to be multiplicative\footnote{\label{foot:multiplicative}That is, $\rho^*(cW)=c\rho^*(W)$ for $c\ge 0$.}, Theorem~\ref{thm:main} tells us that for every $\eps\in(0,1)$, $\OriginalTwoSAT(n,(1-\eps)W)$ and $\OriginalTwoSAT(n,(1+\eps)W)$ are asymptotically almost surely satisfiable and asymptotically almost surely unsatisfiable, respectively. An open question concerning the asymptotic satisfiability of $\OriginalTwoSAT(n,W)$ is given in Section~\ref{ssec:QuestionRhoOne}.
\end{rem}

We complete this section with two propositions which may be useful to combine with Theorem~\ref{thm:main}. Proposition~\ref{prop:contradictorycomponentsproduct} asserts that each strong component $\Omega_i$ of $\overrightarrow{W}$ which is contradictory is in fact as contradictory as it can be. Namely, we have that $\Omega_i=Z\times \signum$ for some $Z\subseteq \Lambda$. This is a counterpart to a well-known and easy property of the implication digraph of a 2-CNF formula\footnote{2-CNF is an abbreviation for 2-Conjunctive Normal Form. These are boolean formula in conjunctive normal form with two literals per clause. So, each instance of random 2-SAT is a 2-CNF formula.} (see Section~\ref{sec:2SATUsingImplicationDigraph} for definition), namely that for each strong component of the implication digraph, we have that it either contains no pair of complementary literals or every literal comes with its complement. 

Proposition~\ref{prop:alternativeexpression} gives an alternative expression for $\rho^*(W)$ by considering all contradictory components combined rather than separately.
\begin{prop}\label{prop:contradictorycomponentsproduct}
Let $W$ be an $L^1$-graphon on $\OurSpace=\Lambda\times\signum$. Let $\OurSpace=\Omega_0\sqcup \bigsqcup_{i\in I}\Omega_i$ be the decomposition of $\overrightarrow{W}$ into strong components. Let $I^*\subseteq I$ be the indices of the contradictory components $\Omega_i$. Then for each $i\in I^*$ we have that $\Omega_i=_0 Z_i\times \signum$ for some $Z_i\subseteq \Lambda$.
\end{prop}
Proposition~\ref{prop:contradictorycomponentsproduct} is proven in Section~\ref{ssec:Proof_contradictorycomponentsproduct}.
\begin{prop}[Proposition~\ref{DIGRAPHONS-prop:spectralradiusAndStrongComponents} and Remark~\ref{DIGRAPHONS-rem:spectralradiusAndStrongComponentsUnbounded} in~\cite{HladkySavicky:Digraphons}]
\label{prop:alternativeexpression}
In the setting of Theorem~\ref{thm:main}, we have 
\[
\sup_{i\in I^*}\rho\left(\overrightarrow{W}\llbracket \Omega_i\rrbracket\right)
=
\max_{i\in I^*}\rho\left(\overrightarrow{W}\llbracket \Omega_i\rrbracket\right)
=\rho\left(\overrightarrow{W}\left\llbracket \bigcup_{i\in I^*}\Omega_i\right\rrbracket\right)=\rho\left(\overrightarrow{W}\left\llbracket \Omega_0\cup \bigcup_{i\in I^*}\Omega_i\right\rrbracket\right)\;.
\]
Moreover, the spectral radii $\rho\left(\overrightarrow{W}\llbracket \Omega_i\rrbracket\right)$ and $\rho\left(\overrightarrow{W}\llbracket \bigcup_{i\in I^*}\Omega_i\rrbracket\right)$ are realized by nonnegative real eigenvalues, so the reference to the modulus in the definition of the spectral radius (Definition~\ref{def:spectrum}) is unnecessary in this case.
\end{prop}

\subsubsection{Theorem~\ref{thm:main} and scale-free random formulea}\label{ssec:scalefree}
Starting in the random graph community around 1999, researchers initiated study of random models in which vertices have heavy-tailed degree distributions. The most famous is arguably the dynamical Barabási--Albert model,~\cite{MR2091634}. For us, the static models of Norros--Reittu~\cite{MR2213964}, and Chung--Lu~\cite{MR2076728} are more relevant. These models became known as `scale-free models', and seem to capture many real-world scenarios. Section~16.4 of~\cite{bollobas2007PhaseTransition} shows that these models are roughly equivalent to the main model $\mathbb{H}(n,U)$ studied~\cite{bollobas2007PhaseTransition} when $U$ is rank-1. We describe $\mathbb{H}(n,U)$ in Section~\ref{ssec:similaritiesToBJR}. Recently, the scale-free phenomenon was introduced also to the setting of random $k$-SAT. Papers~\cite{friedrich2017bounds,DBLP:conf/aaai/0001KRS17} deal with scale-free random $k$-SAT for general $k$, and~\cite{a15060219} deals with 2-SAT. Our Theorem~\ref{thm:main} allows to study a number of scale-free models of random 2-SAT. Translating the main idea from Section~16.4 of~\cite{bollobas2007PhaseTransition} to the setting of random 2-SAT, we could take $\Lambda=(0,1)$, and use kernels of the form $x^{-a}y^{-a}$ (for $a\in (0,1)$) as building bricks. To give a particular example, take $W\in L^1(\OurSpace^2)$ defined as
\begin{align*}
W\big((x,+),(y,+)\big)&=x^{-\alpha}y^{-\alpha}\;, \\
W\big((x,+),(y,-)\big)&=W\big((y,-),(x,+)\big)=x^{-\gamma}y^{-\delta}\;\mbox{, and} \\
W\big((x,-),(y,-)\big)&=(1-x)^{-\beta}(1-y)^{-\beta}\;.
\end{align*}
The purpose of this specific $L^1$-graphon is that it illustrates the limitations on the exponents in Theorem~\ref{thm:main} (which we discussed in Remark~\ref{rem:OnIntegrabilityConditions}). 

We must have $\alpha,\beta,\gamma,\delta<1$ in order to satisfy $W\in L^1(\OurSpace^2)$. As for the additional requirement in Theorem~\ref{thm:main}\ref{en:MainSatisfiable}, we must have $\alpha,\beta,\gamma,\delta<\frac13$ in order to satisfy $W\in L^3(\OurSpace^2)$. As for the additional requirement in Theorem~\ref{thm:main}\ref{en:MainUnsatisfiable}, if we have $\alpha,\beta,\gamma,\delta<\frac12$, then $\overrightarrow{W}\left\llbracket \bigcup_{i\in I^*}\Omega_i\right\rrbracket\in L^2(\OurSpace^2)$ and thus by Lemma~\ref{lem:LpImpliesBoundedCompact}, we have $\overrightarrow{W}\left\llbracket \bigcup_{i\in I^*}\Omega_i\right\rrbracket\in\CompactOper(L^2(\Omega))$.

\subsubsection{Alternative scalings and graphons with $\rho^*(W)=0$}\label{sssec:altarnitivescalings}
In~\eqref{eq:def2SATjb} we defined the probability of the inclusion of any one given clause to be of order $\frac1n$, resulting in $(1\pm o(1)) n \|W\|_1$ clauses in a typical instance of $\OriginalTwoSAT(n,W)$. Some applications might call for other scalings, that is either for sparser models such as $\OriginalTwoSAT(n,(\log n)^{-1}\cdot W)$ or $\OriginalTwoSAT(n,n^{-0.5}\cdot W)$ or for denser models such as $\OriginalTwoSAT(n,\log n\cdot W)$, $\OriginalTwoSAT(n,\sqrt{n}\cdot W)$, or even $\OriginalTwoSAT(n,n\cdot W)$. Note that $\OriginalTwoSAT(n,f(n)\cdot W)$ defines the probability of the clauses using the scaled graphon $\frac{f(n)}{2n} W$. Theorem~\ref{thm:main} together with multiplicativity from Footnote~\ref{foot:multiplicative} has consequences on these other scalings as well. Indeed, first consider the sparser regime $\OriginalTwoSAT(n,f(n)\cdot W)$, where $f(n)\to0$. Clearly, there is $n_0$ such that $f(n)<1/(2\max\{1,\rho^*(W)\})$ for all $n\ge n_0$. We see that the random formula $\OriginalTwoSAT(n,f(n)\cdot W)$ for $n \ge n_0$ is stochastically dominated as a set of clauses by $\OriginalTwoSAT(n,U)$, where $U:=\frac{W}{2\max\{1,\rho^*(W)\}}$. Since $\rho^*(U)\le \frac12$, Theorem~\ref{thm:main} tells us that $\OriginalTwoSAT(n,U)$ is asymptotically almost surely satisfiable, and thus, so is $\OriginalTwoSAT(n,f(n)\cdot W)$. Similarly, in the regime $f(n)\to \infty$, we see that $\OriginalTwoSAT(n,f(n)\cdot W)$ eventually dominates $\OriginalTwoSAT(n,U)$, where $U:=\frac{2W}{\rho^*(W)}$, and thus is asymptotically almost surely unsatisfiable. This argument, however, does not work when $\rho^*(W)=0$. The next proposition resolves this case entirely.
\begin{prop}\label{prop:rhoNula}
Suppose that $W$ is an $L^2$-graphon. The following statements are equivalent.
\begin{enumerate}[label=(\roman*)]
\item\label{en:Dense1} We have $\rho^*(W)=0$.
\item\label{en:Dense2} The $L^2$-digraphon $\overrightarrow{W}$ contains no strong contradictory component.
\item\label{en:Dense3} There exists a function $f:\N\to \R_+$ with $f(n)\to \infty$ so that $\OriginalTwoSAT(n,f(n)\cdot W)$ is not asymptotically almost surely unsatisfiable.
\item\label{en:Dense4} For every function $f:\N\to \R_+$ with $f(n)\to \infty$ we have that $\OriginalTwoSAT(n,f(n)\cdot W)$ is asymptotically almost surely satisfiable.
\end{enumerate}
\end{prop}
As a concrete --- and somewhat surprising --- example, if $\OriginalTwoSAT(n,\log(\log n) \cdot W)$ is asymptotically almost surely satisfiable, then even the considerably denser $\OriginalTwoSAT(n, n \cdot W)$ is as well. We give a proof of Proposition~\ref{prop:rhoNula} in Section~\ref{sec:DenserRegimes}.

\section{Solving 2-SAT using the implication digraph}\label{sec:2SATUsingImplicationDigraph}
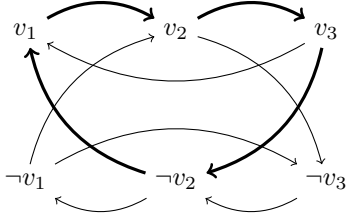
\begin{figure}\centering
\begin{tikzpicture}[->, node distance=2cm, auto]
    \node (x1) {$v_1$};
    \node (x2) [right of=x1] {$v_2$};
    \node (x3) [right of=x2] {$v_3$};

    \node (nx1) [below of=x1] {$\neg v_1$};
    \node (nx2) [below of=x2] {$\neg v_2$};
    \node (nx3) [below of=x3] {$\neg v_3$};

    \path (x1) edge [bend left,very thick] (x2);
    \path (x2) edge [bend left,very thick] (x3);
    \path (x3) edge [bend left] (x1);
    
    \path (nx1) edge [bend left] (x2);
    \path (nx2) edge [bend left,very thick] (x1);
    
    \path (nx3) edge [bend left] (nx2);
    \path (nx2) edge [bend left] (nx1);
    \path (nx1) edge [bend left] (nx3);

    \path (x3) edge [bend left,very thick] (nx2); 
    \path (x2) edge [bend left] (nx3); 
\end{tikzpicture}
\caption{A formula $\phi=(\neg v_1 \vee v_2) \wedge (\neg v_2 \vee v_3) \wedge (\neg v_3 \vee v_1) \wedge (v_1 \vee v_2)\wedge (\neg v_3 \vee \neg v_2)$, and its implication digraph $\ImplDig(\phi)$. One contradictory cycle is highlighted.}
\label{fig:implicationdigraph}
\end{figure}

The concept of the implication digraph of a 2-CNF formula translates the question
of satisfiability into the language of graph theory. Suppose that $\phi$ is
a 2-CNF formula on variables $\{v_1,\ldots,v_n\}$. We create a digraph
\index{$\ImplDig(\phi)$}$\ImplDig(\phi)$, called the \emph{implication digraph of $\phi$}, on the literals $\Lit_n$ as vertices as follows.
For each clause $\{l_1, l_2\}$ we insert directed edges $(\neg l_1, l_2)$
and $(\neg l_2, l_1)$ representing the logically equivalent implications
$\neg l_1 \to l_2$ and $\neg l_2 \to l_1$.
As the edges are inserted in pairs, we see that for each pair $l_i,l_j$ of literals on different variables, $\ImplDig(\phi)$ contains either both directed edges $(l_i, l_j)$ and $(\neg l_j, \neg l_i)$ or none of them. See Figure~\ref{fig:implicationdigraph} for an example.

A directed closed walk or cycle\footnote{in the context of directed graphs,
we consider every
closed walk or cycle directed} $C$ in a digraph $D$ on vertex set $\Lit_n$ is
\emph{contradictory} if there is $i\in[n]$ such that $v_i, \neg v_i \in V(C)$.
It is well-known and easy to verify that a 2-CNF formula $\phi$ is satisfiable if and only if
$\ImplDig(\phi)$ contains no contradictory closed walk. We will use the following stronger statement.
\begin{prop}[\cite{goerdt1996threshold}]\label{prop:SATcycle}
    Suppose that $\phi$ is a 2-CNF formula. Then $\phi$ is satisfiable if and
    only if $\ImplDig(\phi)$ contains no contradictory cycle.
\end{prop}
Proposition~\ref{prop:SATcycle} offers a plausible approach to proving Theorem~\ref{thm:main}. Namely, we need to prove that the implication digraph $\ImplDig(\OriginalTwoSAT(n,W))$ lacks or contains a contradictory cycle in the two respective regimes of Theorem~\ref{thm:main}, asymptotically almost surely. This idea drives the proofs in~\cite{ChvatalR92,goerdt1996threshold} that the satisfiability threshold is $\alpha_2=1$. While the execution works straightforward for some parts of the proof, it turns out that the fact that a contradictory cycle may contain several pairs of complementary literals causes complications in others. To overcome this, \cite{ChvatalR92} introduces a notion of a bicycle. We use a slightly modified form of a bicycle as follows. Suppose that we have integers $k,a,b$, $k \ge 2$, $2 \le a \le k$, and $1 \le b \le k-1$.
Suppose that $D$ is a digraph on $\Lit_n$. A sequence $u_1, \ldots, u_k$ of vertices of $D$ is a \emph{basis of a $(k,a,b)$-bicycle} if the literals $u_1, \ldots, u_k$ are on pairwise
different variables, $u_1 \cdots u_k$ is a directed path
from $u_1$ to $u_k$, and further there is an edge from $\neg u_a$ to $u_1$ and
an edge from $u_k$ to $\neg u_b$ in $D$. We have the following.
\begin{prop}\label{prop:SATbicycle}
    Suppose that $\phi$ is an unsatisfiable 2-CNF formula. Then $\ImplDig(\phi)$
    contains a bicycle which is a subgraph of a contradictory cycle.
\end{prop}
\begin{proof}
By Proposition~\ref{prop:SATcycle}, $\ImplDig(\phi)$ contains a contradictory cycle $C$. Let $P\subseteq C$ be a maximum path in $C$ containing literals on pairwise different variables. $P$ cannot
contain all vertices of the cycle, since $C$ contains complementary literals. $P$ will be the basis of our bicycle. Then, extend $P$ by one edge at each end-point. The result is a bicycle, since $C$ contains no literal twice.
\end{proof}

\subsection{A convenient modification of $\OriginalTwoSAT(n,W)$}\label{ssec:convenientmodel}
We introduce a convenient modification of the model $\OriginalTwoSAT(n,W)$. Here, we recall that $W$ is an  $L^1$-graphon on $\OurSpace$. First, sample $\phi\in\OriginalTwoSAT(n,W)$. We create a formula $\phi^\dagger\sim\TwoSAT(n,W)$ by replacing some of the variables in $\phi$ by their negation. Each variable is replaced by its negation with probability $1/2$ and independently from all other variables.
Note that $\phi$ and $\phi^\dagger$ are equisatisfiable, which means that
either both are satisfiable or both are unsatisfiable. Indeed, we can take a satisfying assignment (if it exists) to one of these formulae, negate the values of the flipped variables, keep the values of the remaining variables, and we get a satisfying assignment of the other formula. The purpose of transforming $\OriginalTwoSAT(n,W)$ to $\TwoSAT(n,W)$ is to simplify analysis of the random formula. The transformation introduces more symmetry into the distribution without affecting the satisfiability we are interested in.

The construction of $\TwoSAT(n,W)$ can be done in a single step as follows,
and we assume this single-step construction in the rest of the paper.
We sample elements $x_1,\ldots,x_n \in \OurSpace$ independently with distribution $\OurMeasure$ and define a map $\tau:\Lit_n \to \OurSpace$ so that for every $i$,
$\tau(v_i)=x_i$ and $\tau(\neg v_i)=\neg \tau(v_i)$. Then, for
every two literals $l_1, l_2 \in \Lit_n$ on different variables, the
clause $\{l_1, l_2\}$ is included in $\TwoSAT(n,W)$ with probability
\begin{equation}\label{eq:2SAT-using-tau}
\min\left\{1,\frac{W(\tau(l_1),\tau(l_2))}{2n}\right\} \;.
\end{equation}
It is easy to verify that the two presented constructions of $\TwoSAT(n,W)$
are equivalent.

\subsection{Random digraph $\G(n,\overrightarrow{W})$}\label{ssec:RandomDigraph}
By Proposition~\ref{prop:SATcycle} we want to get a lower and an upper bound on the probability that the random digraph $\ImplDig(\TwoSAT(n,W))$ contains a contradictory
cycle. By construction, the edges of $\ImplDig(\TwoSAT(n,W))$ are not independent, since
they are included in pairs of the form $(l_1, l_2)$ and $(\neg l_2, \neg l_1)$.\footnote{Recall the definition of the implication digraph in the first paragraph of Section~\ref{sec:2SATUsingImplicationDigraph}.}
In order to overcome
the complications caused by these dependencies in the proof of Theorem \ref{thm:main}\ref{en:MainSatisfiable} we introduce the following random digraph
model \index{$\G(n,U)$}$\G(n,U)$ where $U$ is an $L^1$-digraphon. The vertex set
is $\Lit_n$. Generate a random map $\tau:\Lit_n \to \Omega$ as in
Section~\ref{ssec:convenientmodel} and for each ordered pair of literals
$(l_1, l_2) \in \Lit_n^2$ on different variables insert a directed edge
$(l_1, l_2)$ with probability
\begin{equation} \label{eq:DigraphGeneration}
\min\left\{1,\frac{U(\tau(l_1),\tau(l_2))}{2n}\right\}\;,
\end{equation}
independently of other choices.

Although the distributions of
$\G(n,\overrightarrow{W})$ and of $\ImplDig(\TwoSAT(n,W))$
are closely related, they are not the same.
Indeed, $\G(n,\overrightarrow{W})$ can with positive probability contain a
single edge, whereas $\ImplDig(\TwoSAT(n,W))$ always has an even number of directed edges as explained above. However, restrictions of both these random digraphs on edge sets which do not contain pairs of equivalent edges yield the same marginal distributions. This is true, in particular, for bicycles.
\begin{lem}\label{lem:marginalssame}
Suppose that $W$ is a graphon and $n\in \N$. Let $F\subseteq \Lit_n^2$
be a set with the property that for every pair of literals
$(l_1,l_2) \in \Lit_n^2$ on different variables, we have
$|F\cap \{(l_1, l_2),(\neg l_2, \neg l_1)\}|\le 1$.
Then $\G(n,\overrightarrow{W})_{\restriction F}$ and
$\ImplDig(\TwoSAT(n,W))_{\restriction F}$ have the same distribution. That is, for
every $F'\subseteq F$, we have
\[
\Probability\left[\G(n,\overrightarrow{W})\cap F=F'\right]=\Probability[\ImplDig(\TwoSAT(n,W))\cap F=F']\;.
\]
\end{lem}
\begin{proof}
Recall that both models $\G(n,\overrightarrow{W})$ and $\TwoSAT(n,W)$ start with a random map $\tau:\Lit_n \to \OurSpace$. Let $(l_1,l_2)\in F$ be arbitrary. The probability $\min\left\{1,\frac{\overrightarrow{W}(\tau(l_1),\tau(l_2))}{2n}\right\}$ in~\eqref{eq:DigraphGeneration} for inclusion of the directed edge $(l_1,l_2)$ is the same as the probability $\min\left\{1,\frac{W(\tau(\neg l_1),\tau(l_2))}{2n}\right\}$ for the inclusion of the clause $\{\neg l_1,l_2\}$ in~\eqref{eq:2SAT-using-tau}. This is because of the way $\overrightarrow{W}$ is defined in Definition~\ref{def:implicationdigraphon}. The inclusion of the clause $\{\neg l_1,l_2\}$ is sufficient and necessary for the appearance of directed edge $(l_1,l_2)$ in $\ImplDig(\TwoSAT(n,W))$. 

We conclude that $(l_1,l_2)$ appears as a directed edge in $\G(n,\overrightarrow{W})$ and in $\ImplDig(\TwoSAT(n,W))$ with the same probability. The occurrences of all possible edges in $\G(n,\overrightarrow{W})$, in particular, the edges of $F$, are independent conditionally on $\tau$.
Since the  occurrences of edges of $F$ in $\ImplDig(\TwoSAT(n,W))$ are also independent conditionally on $\tau$, the claim follows.
\end{proof}

\subsection{Basic idea behind the proof of the main theorem}\label{ssec:Heuristics}
We present the main idea behind the proof of Theorem~\ref{thm:main}. We rely on Proposition~\ref{prop:SATcycle}. That is, depending on whether $\rho^*(W)<1$ or $\rho^*(W)>1$, we want to argue that the implication digraph of $\TwoSAT(n,W)$ asymptotically almost surely either does not or does contain a contradictory cycle, respectively. Since any contradictory cycle $C$ itself is a strongly connected directed graph, we have (almost surely) that if $C\subset \ImplDig(\TwoSAT(n,W))$ then all the vertices of $C$ were sampled from one strong component $\Omega_i$ of $\overrightarrow{W}$, that is, (using the notation from Section~\ref{ssec:convenientmodel}) $\tau(l)\in \Omega_i$ for every $l\in V(C)$. Further, since $C$ itself contains a pair of complementary literals, we have that $\Omega_i$ is (almost surely) a contradictory set.

We have therefore reduced the question of asymptotic almost sure nonexistence or existence of a contradictory cycle in $\G(n,\overrightarrow{W})$ to the same question in $\G(n,\overrightarrow{W}\llbracket \Omega_i\rrbracket)$ for a single contradictory strong component $\Omega_i$. It is our task to prove that a contradictory cycle asymptotically almost surely does not exist if $\rho(\overrightarrow{W}\llbracket \Omega_i\rrbracket)<1$ and that it does exist if $\rho(\overrightarrow{W}\llbracket \Omega_i\rrbracket)>1$. For this reason, set $\Gamma:=\overrightarrow{W}\llbracket \Omega_i\rrbracket$ and $\rho:=\rho(\Gamma)$.

Take $a\in \N$ and $b\in \N$. It turns out, that the regime when $a,b=\Theta (\log n)$ is the most relevant one. We first want to get bounds on the probability $p$ that 
\begin{equation}\label{eq:owtt}
v_n,v_1,v_2,\cdots, v_{a-1},v_a,\neg v_n,v_{a+1},v_{a+2},\cdots,v_{a+b}    
\end{equation}
forms a contradictory cycle. We have (ignoring that~\eqref{eq:DigraphGeneration} contains the term $\min\{1,\cdot\}$) that
\[
p=\int_{x_n}\int_{x_1}\ldots\int_{x_{a+b}} \frac{\Gamma(x_{n},x_1)}{2n}\cdot\prod_{s=2}^{a} \frac{\Gamma(x_{s-1},x_s)}{2n}\cdot  \frac{\Gamma(x_{a},\neg x_1)}{2n}\cdot \frac{\Gamma(\neg x_1,x_{a+1})}{2n}\cdot\prod_{s=a+2}^{a+b} \frac{\Gamma(x_{s-1},x_s)}{2n}\cdot \frac{\Gamma(x_{a+b}, x_1)}{2n}\diff\OurMeasure^{1+a+b}\;.
\]
Such products correspond to powers of the digraphon $\Gamma$ (as defined in Definition~\ref{def:digraphonpower}, very much in analogy with powers of square matrices). That is, we have
\[
p=(2n)^{-(a+b+2)}\int_{x_n} \Gamma^{a+1}(x_n,\neg x_n)\cdot \Gamma^{b+1}(\neg x_n,x_n)\diff\OurMeasure\;.
\]
At this point, we use spectral theory. Since a reader might not be familiar with spectral theory of operators on a Hilbert space (and we certainly were not when we started this project), we rather draw an analogy with finite matrices.\footnote{We introduce spectral theory of operators in Section~\ref{ssec:SpectralTheory}.} So, if we could treat $\Gamma$ as a finite matrix with nonnegative entries, the Perron--Frobenius theorem would apply. It says that (under the assumptions of irreducibility and aperiodicity, the former one corresponds to $\Omega_i$ being a strong component, and the latter one is swept under the rug in this sketch) we have $\Gamma^{a+1}(x_n,\neg x_n)=\Theta(\rho^{a+1})$ and $\Gamma^{b+1}(\neg x_n,x_n)=\Theta(\rho^{b+1})$. Note this sketched usage of the Perron--Frobenius theorem relies on the fact that $\Omega_i$ is contradictory. That is, the Perron--Frobenius theorem can only be used in the strongly connected setting, meaning, when there are paths in both directions between $x_n$ and $\neg x_n$. The fact that $\Omega_i$ is contradictory provides this setting for a positive measure of $x_n$'s.

To summarize, we have
\[
p=\Theta(1)\cdot (2n)^{-(a+b+2)}\rho^{a+b+2}\;.
\]

The number of rooted contradictory cycles of length $a+b+2$ with the property that the only pair of complementary literals is the pair 1 and $(a+2)$-nd is equal to $2^{a+b+1}(n)_{a+b+1}$, where $(n)_{a+b+1}$ is a falling factorial. Indeed, we can think of this by replacing each of $v_n,v_1,v_2,\ldots,v_{a+b}$ in~\eqref{eq:owtt} by arbitrary but distinct variables and additionally equipping that variable either with the positive or with the negative sign. Hence, the expected number of such all cycles is
\begin{equation}\label{eq:heuexp}
p\cdot 2^{a+b+1}(n)_{a+b+1}= \Theta(1)\cdot (2n)^{-(a+b+2)}\rho^{a+b+2}\cdot 2^{a+b+1}(n)_{a+b+1}=\Theta(1)\cdot \frac{1}{n}\cdot \rho^{a+b+2}\;.    
\end{equation}
When $\rho<1$ then this quantity is $o(1)$, even when summed over all $a$ and $b$. We conclude that with high probability there are no contradictory cycles. (The above calculation is somewhat simplified since some contradictory cycle might use a several pairs of complementary literals. This is not captured in this sketch.)

When $\rho>1$, then by choosing $a,b=\Theta (\log n)$ large enough, we have $\rho^{a+b+2}\gg n$. Hence the expectation in~\eqref{eq:heuexp} goes to infinity. A second moment argument would be used to show that $\G(n,\Gamma)$ contains a contradictory cycle asymptotically almost surely.

In the actual proof, it is not contradictory cycles that are counted, but different structures: these are bicycles (for the proof of Theorem~\ref{thm:main}\ref{en:MainSatisfiable}) and what is later called `snakes' (for the proof of Theorem~\ref{thm:main}\ref{en:MainUnsatisfiable}). 

\subsection{Comparison to previous work on inhomogeneous random (di)graphs}\label{ssec:similaritiesToBJR}

As discussed in Section~\ref{ssec:Heuristics}, a central aspect of our proof lies in estimating the probability that certain special (`contradictory') cycles occur in inhomogeneous random digraphs $\G(n,\overrightarrow{W})$. In this section, we review several foundational results on inhomogeneous random graphs and digraphs, which concern the appearance of a giant component. We will remark at the end of this section that, while our main results concern contradictory cycles, these two phenomena are likely related, though we do not rely on this connection in our arguments.

A cornerstone of random graph theory is the well-known phase transition in the Erdős–Rényi random graph on $n$ vertices with edge probability $\frac{c}{n}$. This model exhibits a sharp threshold at $c=1$: when $c<1$, the largest connected component is, with high probability, of order $o(n)$, whereas for $c>1$, a unique giant component of order $\Theta(n)$ emerges.

Bollobás, Janson, and Riordan~\cite{bollobas2007PhaseTransition} extended this phenomenon to a sparse inhomogeneous random graph model $\mathbb{H}(n,U)$, parameterized by an almost everywhere continuous $L^1$-graphon $U$ defined on a metric probability space $\Lambda$ with measure $\lambda$. In their simplest setting, the vertex set of $\mathbb{H}(n,U)$ is $[n]$, and vertices $x_1,\dots,x_n \in \Lambda$ are sampled independently according to~$\lambda$. Each edge $ij$ is then included independently with probability $\min\{1,U(x_i,x_j)/n\}$. This model and its variants have since been widely studied. The general framework of~\cite{bollobas2007PhaseTransition} also allows for more flexible vertex-type sampling and for random perturbations depending on~$n$, but such generality requires $U$ to be almost everywhere continuous.

Our random digraph model $\G(n,\overrightarrow{W})$, introduced in Section~\ref{ssec:RandomDigraph}, is closely related to this framework, or more precisely, to its directed analogue. This directed setting was first analyzed in a restricted case by Bloznelis, Götze, and Jaworski~\cite{ND}, and later extended by Cao and Olvera-Cravioto~\cite{InhomogeneousRandomDigraphs}. We next summarize their results and explain their relevance to our work.

A key theorem of~\cite{bollobas2007PhaseTransition} provides a criterion for the existence of a giant component in $\mathbb{H}(n,U)$. Specifically, if $\rho(U)\le 1$, then with high probability $\mathbb{H}(n,U)$ has no giant component as $n\to\infty$; conversely, if $\rho(U)>1$, a giant component appears. The proof proceeds in two steps: first, the local neighborhood of a vertex in $\mathbb{H}(n,U)$ is shown to correspond to an inhomogeneous Galton–Watson branching process~$\mathcal{S}_U$; second, the condition $\rho(U)\le 1$ implies almost sure extinction of~$\mathcal{S}_U$, whereas $\rho(U)>1$ yields a positive survival probability. These arguments, presented in Section~5 of~\cite{bollobas2007PhaseTransition}, are self-contained.

Bloznelis, Götze, and Jaworski~\cite{ND} developed an analogous branching-process correspondence for inhomogeneous random digraphs, restricted to stochastic block models, though without invoking spectral arguments. Cao and Olvera-Cravioto~\cite{InhomogeneousRandomDigraphs} extended the analysis to a broader family of directed models $\G(n,\Gamma)$, allowing for different type samplings and random perturbations in the spirit of~\cite{bollobas2007PhaseTransition}. Under the assumption that $\Gamma$ is almost everywhere continuous, their analysis effectively reduces to the stochastic block model setting, while the functional-analytic component of their argument, as in~\cite{MR4011679}, remains relatively elementary.

Finally, we note a heuristic link between cycles and the giant component phenomenon. In the Erdős–Rényi model, it is well known that sublinear components contain only few cycles, whereas the giant component contains exponentially many.
Thus, the appearance of a giant component and the appearance of cycles
through specific types of vertices are likely connected, although our proof
does not rely on this connection.
Besides this, we introduce advanced spectral techniques that
allow us to treat general (non-continuous) digraphons $\Gamma$ directly, without relying on stochastic block model approximations.

\section{Preliminaries}
\subsection{Measure theory}
While parts of the paper specific to 2-SAT require to work in a the product space $\OurSpace=\Lambda\times \signum$, it is convenient to formulate some more basic parts with respect to a general probability space. We use a measure space \index{$(\Omega,\mu)$}$\Omega$ equipped with a probability measure $\mu$ to this end. That is, the measure $\mu$ always implicitly underlies the space $\Omega$.

For measurable subsets $A$ and $B$ in a measure space, we write \index{$=_0$, $\subseteq_0$}$A=_0 B$ and $A\subseteq_0 B$ for equality and containment modulo a nullset, respectively. We write \index{$\essinf f$, $\esssup f$}$\essinf f$ and $\esssup f$ for the essential infimum and essential supremum of a function $f$ on a measure space.

\subsection{Digraphons}\label{ssec:Digraphons}
Suppose that $p\in [1,\infty)$, $\Gamma$ is an $L^p$-digraphon on $\Omega$ and $D$ is an oriented graph on vertex set $[n]$. The \emph{homomorphism density of $D$ in $\Gamma$} is defined as
\begin{equation}\label{eq:defHomDens} 
\index{$t(D,\Gamma)$}
t(D,\Gamma)=\int_{x_1}\int_{x_2}\ldots\int_{x_n}\prod_{(i,j) \in E(D)}\Gamma(x_i,x_j)\diff\mu^V\;.
\end{equation}
This quantity may be infinite, but is certainly finite when $\Gamma$ is a digraphon.

Suppose that $k\in \N$. Write \index{$P_k,P_k^{\bullet\bullet}$}$P_k$ for the directed path $1,2,3,\cdots ,k,{k+1}$, and $P_k^{\bullet\bullet}$ for the directed path rooted at its terminal vertices. Write \index{$C_k,C_k^\bullet$} for the directed cycle $1,2,3,\cdots ,k$ and $C_k^\bullet$ for $C_k$ rooted at vertex $1$. We define \index{$t^{\bullet\bullet}_{x,y}(P_k^{\bullet\bullet},\Gamma)$}$t^{\bullet\bullet}_{x_1,x_{k+1}}(P_k^{\bullet\bullet},\Gamma)$ as a function of $x_1,x_{k+1}\in\Omega$ by disintegrating $t(P_k,\Gamma)$ with respect to $x_1$ and $x_{k+1}$. Likewise, we define \index{$t^{\bullet}_{x}(C_k^{\bullet},\Gamma)$}$t^{\bullet}_{x_1}(C_k^{\bullet},\Gamma)$ as a function of $x_1\in\Omega$ by disintegrating $t(C_k,\Gamma)$ with respect to $x_1$. That is,
\begin{align}    
\label{eq:densPathTwoTerminals}
t^{\bullet\bullet}_{x_1,x_{k+1}}(P_k^{\bullet\bullet},\Gamma)&=\int_{x_2}\int_{x_3}\ldots\int_{x_k}\prod_{i=1}^{k}\Gamma(x_i,x_{i+1})\;\diff\mu^{k-1}\;,\\
\nonumber
t^{\bullet}_{x_1}(C_k^{\bullet},\Gamma)&=\int_{x_2}\int_{x_3}\ldots\int_{x_k}\prod_{i=1}^{k-1}\Gamma(x_i,x_{i+1})\cdot \Gamma(x_k,x_1)\;\diff\mu^{k-1}\;.
\end{align}

Next, we introduce the power of a digraphon. This definition is similar to the definition of matrix powers, and is also used in connection with integral kernel operators.
\begin{defi}[Power of a digraphon]\label{def:digraphonpower}
Suppose that $\Gamma$ is a digraphon on $\Omega$, and let $k\in\N$. Define \index{$\Gamma^k$}$\Gamma^k$ as digraphon on $\Omega$ by $\Gamma^k(x,y):=t^{\bullet\bullet}_{x,y}(P_k^{\bullet\bullet},\Gamma)$.
\end{defi}
It is straightforward to verify that the definition is consistent with operator powers, that is, for every $f\in L^2(\Omega)$ and every $k\in \N$ we have
\[
\Gamma^k f=\underbrace{\Gamma(\Gamma(\cdots (\Gamma}_{\mbox{$k$ times}} (f))\cdots))\;.
\]
\begin{rem}
We will frequently use the concept of a power in the sense of Definition \ref{def:digraphonpower}, as well as the notion of taking a power of the value of a digraphon. To distinguish between them, we write $\Gamma^k(x,y)$ for the former and $\Gamma(x,y)^k$ for the latter. For example, if $\Gamma$ is a digraphon defined on the unit square, $\Gamma(x,y):=\mathbbm{1}_{y\le 0.5}$, then $\Gamma^k(0.1,0.2)=0.5^{k-1}$ and $\Gamma(0.1,0.2)^k=1$ for all $k\in \N$.
\end{rem}

\subsection{Banach space theory}\label{ssec:SpectralTheory}
Our proof of Theorem~\ref{thm:main} uses a fair amount of functional analysis. We summarize the tools we need in a way which should be accessible to discrete mathematicians.

As we said earlier, all our Hilbert/Banach spaces will be complex. We write \index{$\ImaginaryUnit$}$\ImaginaryUnit$ for the imaginary unit. Specifically, we will work with Banach spaces $L^p(\Omega)$ (for $p\in [1,\infty)$) and Hilbert space $L^2(\Omega)$. For a Banach space $X$, we write \index{$\BoundOper(X)$}$\BoundOper(X)$ and \index{$\CompactOper(X)$}$\CompactOper(X)$ for the set of bounded operators $X\to X$ and the set of compact operators $X\to X$, respectively.

\begin{defi}\label{def:otherLp}
Suppose that $p\in[1,\infty)$ is given. For an operator $T\in\BoundOper(L^{p}(\Omega))$,
define its \emph{operator norm} \index{$\|T\|_{\mathrm{op}(p)}$}$\|T\|_{\mathrm{op}(p)}=\sup_{f\in L^{p}(\Omega),\|f\|_{p}=1}\|Tf\|_{p}$.
\end{defi}
Observe that for every $f\in L^{p}(\Omega)$, 
\begin{equation}\label{eq:SJI}
\|Tf\|_{p}^{p}\le\|T\|_{\mathrm{op}(p)}^{p}\cdot\|f\|_{p}^{p}\;.
\end{equation}

In Definition~\ref{def:spectrum}, we defined eigenvalues and eigenfunctions of an $L^2$-digraphon $\Gamma$. These eigenvalues and eignefunctions are sometimes also called \emph{right}. \emph{Left eigenvalues} and \emph{left eigenfunctions} of $\Gamma$ are (right) eigenvalues and eigenfunctions of transposed $L^2$-digraphon $\Gamma^{\top}$, $\Gamma^{\top}(x,y):=\Gamma(y,x)$.

The lemma below is standard and concerns digraphons as integral kernel operators.
\begin{lem}\label{lem:eigenvectorsbounded}
Suppose that $\Gamma$ is a digraphon on $\Omega$. Then for arbitrary $p\in[1,\infty)$, we have $\Gamma\in\BoundOper(L^p(\Omega))$. Further, consider an arbitrary eigenvalue $\gamma\neq 0$ and corresponding eigenfunction $f$, when $\Gamma$ is viewed as an operator on $L^p(\Omega)$. Then $\|f\|_\infty\le \frac{\|\Gamma\|_\infty \|f\|_1}{|\gamma|}$. Consequently, $\PSpec_p(\Gamma)$ and $\rho_p(\Gamma)$ does not depend on the choice of $p\in[1,\infty)$.
\end{lem}
\begin{proof}
First, we show that $\Gamma\in\BoundOper(L^p(\Omega))$. Consider an arbitrary $h\in L^p(\Omega)$. We have $\|\Gamma h\|_p^p=\int_z\int_x |\Gamma(z,x)h(x)|^p\le \|\Gamma\|_\infty^p\int_z|h(z)|^p<\infty$. Hence, $\Gamma h\in L^p(\Omega)$.

We now turn to $L^\infty$-boundedness of the eigenvalues. We have $\gamma f=\Gamma f$. Fix an arbitrary $x\in\Omega$. Thus, $|\gamma f(x)|=|\int_y f(y)\Gamma(y,x)|\le \|\Gamma\|_\infty\cdot\|f\|_1$.
\end{proof}

\subsubsection{\Krein--Rutman theorem and beyond}
Suppose that $\Gamma$ is an digraphon on $\Omega$ and $\tau$ is its eigenvalue. We say that $\tau$ is \emph{simple} if for every $f_1,f_2\in L^2(\Omega)$ with $\Gamma f_i=\tau f_i$ we have that $f_1$ is a scalar multiple of $f_2$.

We now introduce a version of the \Krein--Rutman Theorem. Let us give some background first. The Perron--Frobenius Theorem asserts that a real square matrix with positive entries has a unique nonnegative eigenvector, and that this eigenvector corresponds to the eigenvalue of the largest eigenvalue in absolute value. The \Krein--Rutman Theorem is often considered a counterpart of the Perron--Frobenius Theorem for nonnegative operators. While there are many version, we reproduce a version tailored to our language of digraphons, included as Theorem~\ref{DIGRAPHONS-thm:Schaefer74} in~\cite{HladkySavicky:Digraphons}.
\begin{thm}\label{thm:Schaefer74}
Suppose that $\Gamma$ is a strongly connected digraphon on $\Omega$ with $\rho(\Gamma)>0$.
Then for some $d \in\N$, the set of eigenvalues of maximal modulus is
$\left\{\exp(-2\pi \ImaginaryUnit k/d) \rho(\Gamma)\::\:k=0,\ldots,d-1\right\}$. Further, there are right and left  eigenfunctions $v_\mathrm{right}$ and $v_\mathrm{left}$ for the eigenvalue $\rho(\Gamma)$. Both these eigenfunctions are strictly positive and we have $v_\mathrm{right},v_\mathrm{left}\in L^\infty(\Omega)$.
\end{thm}
We call the number $d$ from Theorem~\ref{thm:Schaefer74} the \emph{peripheral multiplicity} of $\Gamma$.

We use the famous Gelfand's formula  (see P7.5-5 in~\cite{MR0467220}). Recall that the \emph{operator norm} $\|R\|_{\mathrm{op}}$ of an operator $R$ on a Banach space $(X,\|\cdot\|)$ is defined as \index{$\|\cdot \|_{\mathrm{op}}$}$\|R\|_{\mathrm{op}}:=\sup_{x\in X:\|x\|\le 1}\|Rx\|$.
\begin{prop}\label{pro:Gelfand}
Suppose that $T$ is a bounded operator on a complex Banach space. Then we have $\rho(T)=\lim_{k\to\infty}\sqrt[k]{\|T^k\|_{\mathrm{op}}}$.
\end{prop}
Gelfand's formula holds for a variety of other norms. We recall another version, taken from~\cite{HladkySavicky:Digraphons} (though, it is very likely to be known). This version concerns the Hilbert--Schmidt norm, $\|\cdot\|_{\mathrm{HS}}$. As we are concerned with $L^2$-digraphons, we shall need the Hilbert--Schmidt norm only in the case of integral kernel operators, and thus we introduce it only in this particular setting. It is well-known that if $T_K$ is an integral kernel operator on $L^2(\Omega)$ with kernel $K\in L^2(\Omega^2)$ then $\|T_K\|_{\mathrm{HS}}=\|K\|_2$.
\begin{prop}[Proposition~\ref{DIGRAPHONS-pro:GelfandHS} in~\cite{HladkySavicky:Digraphons}]\label{pro:GelfandHS}
Suppose that $\Gamma\in L^2(\Omega^2)$ is an $L^2$-digraphon. Then we have $\rho(\Gamma)=\lim_{k\to\infty}\sqrt[k]{\|\Gamma^k\|_{2}}$.
\end{prop}

Gelfand's formula (Proposition~\ref{pro:Gelfand}) asserts that 
$\rho_{p}(T)=\lim_{k\to \infty}\left(\left\Vert T^{k}\right\Vert _{\mathrm{op}(p)}\right)^{1/k}$.
No additional assumptions on $T$ are needed. 
The following lemma says that the spectral radius of a nonnegative
operator is independent of the Banach space we work with in many scenarios. This is
well known and we recall its proof for completeness. The statement and the proof of the lemma uses some basic notions (nonnegativity, irreducibility) from the theory of Banach lattices. The reader can find these notions explained in our language in Section~\ref{DIGRAPHONS-ssec:BanachLattices} of~\cite{HladkySavicky:Digraphons}. What is important to us is that the nonnegativity of digraphons implies that they are nonnegative on each Banach lattice $L^{q}(\Omega)$.
\begin{lem}
\label{lem:GelfandInvariant}
Suppose that $1\le q\le p<\infty$ are two numbers and $T$ is a bounded operator, both as $T:L^{q}(\Omega)\to L^{q}(\Omega)$ and as $T:L^{p}(\Omega)\to L^{p}(\Omega)$.
\begin{enumerate}[label=(\roman*)]
\item\label{en:lGI1}
We have $\rho_q(T)\ge \rho_p(T)$.
\item\label{en:lGI2}
If $T$ is nonnegative (in the sense
of Banach lattices) and $T\in \CompactOper(L^p(\Omega))\cap \CompactOper(L^{q}(\Omega))$, then $\rho_{p}(T)=\rho_{q}(T)$.
\end{enumerate}
\end{lem}
\begin{proof}
Part~\ref{en:lGI1} is simple. Each eigenfunction of $T:L^{p}(\Omega)\to L^{p}(\Omega)$ is also an eigenfunction of $T:L^{q}(\Omega)\to L^{q}(\Omega)$ (with the same eigenvalue).

Let us turn to Part~\ref{en:lGI2}. Without loss of generality, let us assume $T$ is irreducible (if it
were not, we would break it into irreducible parts). Let $f$ be the
eigenfunction corresponding to the eigenvalue of the maximum modulus
when $T$ is viewed as an operator on $L^{p}(\Omega)$, and let $g$
be the eigenfunction corresponding to the eigenvalue of the maximum
modulus when $T$ is viewed as an operator on $L^{q}(\Omega)$. We will prove that $f$ and $g$ are the same (up to a multiple), which will prove the statement. The \Krein-Rutman theorem (applied
twice, once for operators on $L^{p}(\Omega)$ and once for operators
on $L^{q}(\Omega)$) tells us that $f$ and $g$ are in the positive
cone (in other words, they are nonnegative functions). Also, there
is a uniqueness part of the \Krein--Rutman theorem, which asserts that
(for irreducible operators), there is only one nonnegative eigenfunction.
But in the space $L^{p}(\Omega)$ we seem to have two nonnegative
eigenfunctions, namely $f$ and $g$. So, this is only possible $f$ and $g$ are the same (up to a multiple).
\end{proof}

\subsubsection{Approximating an integral kernel operator}
In Section~\ref{sec:Unsatisfiable}, we prove Theorem~\ref{thm:main}\ref{en:MainUnsatisfiable}. The calculations there do not work in the general case $W\in L^1(\Omega^2)$. So, as a first step in Section~\ref{sec:Unsatisfiable}, we find a digraphon $V\in L^\infty(\OurSpace^2)$ with $V\le W$ and $\rho^*(V)>1$. This assumption of $L^\infty$-boundedness then allows us to compute that indeed $\OriginalTwoSAT(n,V)$ is asymptotically almost surely unsatisfiable. Since $V\le W$, the same is true for $\OriginalTwoSAT(n,W)$. The following result is used for said approximation from below.
\begin{prop}\label{prop:approximatefrombelow}
Suppose that we have $p\in[1,\infty)$ and $W$ an $L^1$-graphon on $\OurSpace$ with $\overrightarrow{W}\in\CompactOper(L^p(\Omega))$. Then for every $\eps>0$, there exists a graphon $V$ on $\OurSpace$ with the properties that $V \le W$, and $\rho^*(V)\ge \rho^*_p(W)-\eps$.
\end{prop}
\begin{proof}
We define a sequence of graphons $V_1\le V_2\le \ldots \le W$, where $V_n$ is the pointwise minimum of $W$ and the constant-$n$ function, $V_n:=\min(W,n)$. Obviously, $V_n$ is symmetric and bounded, and thus a graphon. Also, the decomposition of $\overrightarrow{V_n}$ into strong components is the same as of $\overrightarrow{W}$. So, the sets of contradictory components involved in the definition of $\rho^*_p(V_n)$ and of $\rho^*_p(W)$ are the same.
Thus, we have $\lim_{n\to\infty} \rho^*_p(V_n)=\rho^*_p(W)$ by Theorem~2.4 in~\cite{Schep} (applied to the Banach space $L^p(\OurSpace)$). We conclude that there exists a graphon $V\le W$ with $\rho^*_p(V)\ge \rho^*_p(W)-\eps$. Since $\overrightarrow{V}$ is a digraphon, Lemma~\ref{lem:eigenvectorsbounded} applies. It tells us that the notion of the spectral radius does not depend on the choice of the space $L^q(\Omega)$ (over all $q\in[1,\infty)$). Therefore, we have $\rho^*_p(V)=\rho^*(V)$.
\end{proof}

\subsubsection{$L^p$-digraphons as integral kernel operators on $L^p(\Omega)$}\label{sssec:GelfandLp}
Suppose that $W$ is an $L^p$-digraphon (in this section, we work with $p\in[2,\infty)$). In this section, we deduce that $W$ as an integral kernel operator on $L^p(\Omega)$ is bounded and compact (Lemma~\ref{lem:LpImpliesBoundedCompact}) and express the spectral radius using a Gelfand-like formula involving the $L^p$-norm. While we were not able to find this result, we believe it could be known.

The first result is similar to Exercise~7 on page~177 in~\cite{MR1070713}.
\begin{lem}\label{lem:LpImpliesBoundedCompact}
Suppose that $W\in L^p(\Omega^2)$ for some $p\ge 2$. Then $W$ as an integral kernel operator on $L^p(\Omega)$ is bounded with $\|W\|_{\mathrm{op}(p)}\le \|W\|_p$. Furthermore, $W$ is compact on $L^p(\Omega)$.
\end{lem}
\begin{proof}
Suppose that $f\in L^p(\Omega)$. First, we prove that $Wf$ is well-defined as a function at almost every $x\in \Omega$. Let $q$ be the Hölder conjugate to $p$, $\frac1p+\frac1q=1$. We have that $q\le 2\le p$. In particular, $\|f\|_q\le \|f\|_p$. Also, notice that when $\|W\|_p<\infty$, then by Fubini's theorem, at almost every $x\in \Omega$ we must have $\|W(x,\cdot)\|_p<\infty$. Hence, using Hölder's inequality (HI),
\[
|(Wf)(x)|=\left|\int W(x,y)f(y)\diff\mu(y)\right|\leBy{(HI)} \|W(x,\cdot)\|_p\cdot\|f\|_q\le \|W(x,\cdot)\|_p\cdot\|f\|_p<\infty\;.
\]
Let us now get a bound on the $L^p$-norm of $Wf$. We use the previous calculation
\[
\|Wf\|^p_p = \int |(Wf)(x)|^p\diff\mu(x) \le \int \|W(x,\cdot)\|_p^p\cdot\|f\|_p^p\diff\mu(x)\eqBy{Fubini}\|W\|_p^p\cdot\|f\|_p^p \;.
\]
Equivalently, $\|Wf\|_p\le \|W\|_p\cdot\|f\|_p$. This shows that $W$ as an integral kernel operator on $L^p(\Omega)$ is bounded with $\|W\|_{\mathrm{op}(p)}\le \|W\|_p$.

We now turn to proving compactness of $W$. Recall that the product sigma-algebra on $\Omega^2$ is generated by sets of the form $S\times T$, where $S,T\subset\Omega$ are measurable. That is, it is known that for an arbitrary $\eps>0$, we can find a finite sequence $S_1,T_1,S_2,T_2,\ldots,S_\ell,T_\ell\subset \Omega$ and coefficients $c_1,c_2,\ldots,c_\ell\in \R$ so that for the function $W_\eps:=\sum_{i=1}^\ell c_i \mathbbm{1}_{S_i\times T_i}$ we have $\|W-W_\eps\|_p<\eps$. By the previous (applied to $U=W-W_\eps$), we have $\|W-W_\eps\|_{\mathrm{op}(p)}\le \|W-W_\eps\|_p<\eps$. Also, the rank of $W_\eps$ is at most $\ell$. That is, up to arbitrary precision, we are able to approximate $W$ by a finite-rank operator in the operator norm. By a well-known fact (see e.g. Theorem~4.4 in~\cite{MR1070713}), it follows that $W$ is compact.
\end{proof}

There are many versions of Gelfand's formula (Proposition~\ref{pro:Gelfand}). We will need a different version, involving $L^p$-norms of kernels of a nonnegative $L^p$-digraphon.
\begin{prop}\label{prop:GelfandLp}
Let $p\in[2,\infty)$. Suppose that $W$ is an $L^p$-digraphon on $\Omega$. Then we have $\rho_p(W)=\lim_{k\to\infty}\left(\left\Vert W^{k}\right\Vert _{p}\right)^{1/k}$.
\end{prop}

The next lemma is an important ingredient for the proof of Proposition~\ref{prop:GelfandLp}. For it, we work with \emph{composition of kernels}. That is, if $A,B\in L^{p}(\Omega^{2})$, then we define $C:=A*B$ as a function on $\Omega^2$, $C(x,y):=\int_z A(x,z)B(z,y)$ (provided that the integral is defined). The next lemma tells us that under mild conditions, $A*B$ is indeed defined almost everywhere, and we have $A*B\in L^{p}(\Omega^{2})$.
\begin{lem}
\label{lem:Jep}
Let $p\in[1,\infty)$ be arbitrary.
Let $A,B\in L^{p}(\Omega^{2})$ be two kernels. Let $T_A$ be the integral kernel operator associated with $A$, and
assume that $T_A\in \BoundOper(L^{p}(\Omega))$. Set $C=A*B$.
Then $C$ is defined almost everywhere, $C\in L^p(\Omega^2)$, and $\left\Vert C\right\Vert _{p}^{p}\le\left\Vert T_{A}\right\Vert _{\mathrm{op}(p)}^{p}\cdot\left\Vert B\right\Vert _{p}^{p}$.
\end{lem}
\begin{proof}
For $x\in\Omega$, write $C_{x}$ for the ``column slice of $C$
at $x$''. That is, $C_{x}$ is a one-variable function, $C_{x}(y):=C(x,y)$.
Fubini's Theorem gives
\[
\left\Vert C\right\Vert _{p}^{p}=\int_{x}\left\Vert C_{x}\right\Vert _{p}^{p}\diff\mu(x)\;.
\]
With a similar slice convention for $B$, we have
\[
\left\Vert B\right\Vert _{p}^{p}=\int_{x}\left\Vert B_{x}\right\Vert _{p}^{p}\diff\mu(x)\;.
\]
We will prove the lemma slice-wise, that is, we will prove that
for every $x\in\Omega$, $\left\Vert C_{x}\right\Vert _{p}^{p}\le\left\Vert T_{A}\right\Vert _{\mathrm{op}(p)}^{p}\cdot\left\Vert B_{x}\right\Vert _{p}^{p}$.
As $C=A*B$, we have $C_{x}=T_{A}(B_{x})$. So, the inequality
we need is just~(\ref{eq:SJI}) applied to the function $B_{x}$ and
operator $T_{A}$.
\end{proof}

\begin{proof}[Proof of Proposition~\ref{prop:GelfandLp}]
To get one inequality of the proposition, we use Gelfand's formula
for the Hilbert--Schmidt norm (Proposition~\ref{pro:GelfandHS}) and then use that $\left\Vert \cdot\right\Vert _{2}\le\left\Vert \cdot\right\Vert _{p}$ and $\rho(\cdot)=\rho_2(\cdot)\ge \rho_p(\cdot)$ (see~Lemma~\ref{lem:GelfandInvariant}\ref{en:lGI1}).
\[
\rho_p(W)\le \rho(W)=\lim_{k\to\infty}\left(\left\Vert W^{k}\right\Vert _{2}\right)^{1/k}\le\liminf_{k\to\infty}\left(\left\Vert W^{k}\right\Vert _{p}\right)^{1/k}.
\]

We now turn to the other inequality. Lemma~\ref{lem:Jep} with $A=W^{k-1}$
and $B=W$ gives
\begin{align*}
\limsup_{k\to\infty}\left(\left\Vert W^{k}\right\Vert _{p}\right)^{1/k} & =\limsup_{k\to\infty}\left(\left\Vert W^{k}\right\Vert _{p}^{p}\right)^{1/(pk)}\le\limsup_{k\to\infty}\left(\left\Vert T_{W}^{k-1}\right\Vert _{\mathrm{op}(p)}^{p}\cdot\left\Vert W\right\Vert _{p}^{p}\right)^{1/(pk)}\\
 & =\limsup_{k\to\infty}\left(\left\Vert T_{W}^{k-1}\right\Vert _{\mathrm{op}(p)}\right)^{1/k}\cdot\lim_{k\to\infty}\left(\left\Vert W\right\Vert _{p}^{p}\right)^{1/(pk)}\;.
\end{align*}
The first term goes to the spectral radius $\rho_{p}(W)$ by Gelfand's formula (Proposition~\ref{pro:Gelfand}).
The second term goes to 1. Hence, $\limsup_{k\to\infty}\left(\left\Vert W^{k}\right\Vert _{p}\right)^{1/k}\le\rho_p(W)\cdot1$,
as was needed.
\end{proof}

\subsection{More on digraphons}
\subsubsection{Peripheral multiplicity and graphic periodicity of digraphons}
In a directed graph, the \emph{period} is defined as the greatest common divisor of the lengths of all its directed cycles.  
It is well-known that when this period is greater than $1$, the vertex set can be partitioned into equivalence classes called \emph{cyclic sets}.  
Edges in the digraph move vertices from one cyclic set to the next in a fixed cyclic order modulo the period.  
If the period is $1$, the digraph is called \emph{aperiodic}, and no nontrivial cyclic partition exists.
Similar notion exists in the theory of Markov chains.
In~\cite{HladkySavicky:Digraphons}, a counterpart for digraphons was introduced.
\begin{defi}\label{def:graphicallyperiodic}
    Suppose that $\Gamma$ is a digraphon on $\Omega$. For $d\in\N$, we say that $\Gamma$ is \emph{graphically $d$-periodic} if there exists a partition $\Omega=X_0\sqcup X_1\sqcup \ldots\sqcup X_{d-1}$ such that (using the cyclic notation $X_d=X_0$) for every $j=0,\ldots,d-1$ we have $\Gamma_{\restriction X_j\times (\Omega\setminus X_{j+1})}=0$.
\end{defi}
The following was proven in~\cite{HladkySavicky:Digraphons}.
\begin{thm}[Theorem~\ref{DIGRAPHONS-thm:periodicity} in~\cite{HladkySavicky:Digraphons}]\label{thm:periodicity}
Suppose that $\Gamma$ is a strongly connected digraphon on $\Omega$. Suppose that the peripheral multiplicity of $\Gamma$ is $D$. 
Then $\Gamma$ is graphically $D$-periodic.
\end{thm}

\subsubsection{Key asymptotics}
The following is one of the main results of~\cite{HladkySavicky:Digraphons}. It asserts that high powers of a digraphon can be asymptotically expressed using the spectral radius, the left principal eigenfunction and the right principal eigenfunction.
\begin{prop}[Theorem~\ref{DIGRAPHONS-thm:asymptotics} in~\cite{HladkySavicky:Digraphons}]
 \label{prop:asymptoticpowers}
Suppose that $\Gamma$ is a strongly connected digraphon on ground set $\Omega$.
We assume that there are left and right real eigenfunctions
$v_L$, $v_R$ for the eigenvalue $\rho(\Gamma)$ satisfying $\langle v_L, v_R \rangle = 1$.

Let the peripheral multiplicity of $\Gamma$ be $D$. Suppose that $\Omega=X_0\sqcup X_1\sqcup \ldots\sqcup X_{D-1}$ is a decomposition as in Definition~\ref{def:graphicallyperiodic} provided by Theorem~\ref{thm:periodicity}.

Let $\rho:=\rho(\Gamma)$. There exists a number $\alpha\in (0,\rho)$ with the following property. For every $i,j\in\{0,\ldots,D-1\}$ and every $x\in X_i$ and $y\in X_j$ we have
\begin{equation*}
\Gamma^\ell(x,y) = \begin{cases}
                    \rho^\ell v_R(x) v_L(y) + O(\alpha^\ell)&\quad\mbox{if $\ell\equiv j-i \mod D$, or}\\
                    0&\quad\mbox{otherwise,}
                   \end{cases} 
\end{equation*}
as $\ell\to\infty$. The term in $O(\cdot )$ does not depend on $x$ and $y$.
\end{prop}

\subsection{Proof of Proposition~\ref{prop:contradictorycomponentsproduct}}\label{ssec:Proof_contradictorycomponentsproduct}
Suppose that $X\subset \OurSpace$ is a strong component of $\overrightarrow{W}$. Since $X$ is a contradictory component, the set $P:=\{z\in \OurSpace:z,\neg z\in X\}$ has positive measure. Define $R:=X\setminus P$. Define $\overline{R}:=\{\neg z:z\in R\}$. We shall prove that 
\begin{equation}\label{eq:dasenka}
\mbox{$P\sqcup R\sqcup \overline{R}$ is strongly connected.}
\end{equation}
From that it will follow, using the fact that $X$ is a strong component and Definition~\ref{def:component}\ref{en:StrongComp}, that $\overline{R}$ is null, which in turn yields that $R$ is null. This will prove the statement.

To prove~\eqref{eq:dasenka}, consider a partition $A\sqcup B=P\cup R\cup \overline{R}$ into two arbitrary sets $A$ and $B$ of positive measure, as in Definition~\ref{def:component}\ref{en:StrongConn}. We have
$
\int_{A\times B}\overrightarrow{W}\ge \int_{(A\cap X)\times (B\cap X)}\overrightarrow{W}
$. The last term is positive by the fact that $X$ is a strong component and Definition~\ref{def:component}\ref{en:StrongConn}, provided that we prove that $A\cap X$ and $B\cap X$ have positive measure. 
So, it only remains to deal with the cases that $A\cap X$ or $B\cap X$ is null. Suppose for example that $A\cap X$ is null, the other case being analogous. Then $A\subset_0 \overline{R}$. In particular, for the set $\overline{A}:=\{\neg z:z\in A\}$ we have $\overline{A}\subset R$. Also, the set $\overline{B}:=\{\neg z:z\in B\}$ contains $P$. That means that $(\overline{B}\cap X)\sqcup \overline{A}$ is a partition of $X$ into two sets of positive measures. By~\eqref{eq:antisymmetric}, we have $\int_{(A\cap X)\times (B\cap X)}\overrightarrow{W}=\int_{(\overline{B}\cap X)\times \overline{A}}\overrightarrow{W}>0$, where the last inequality follows from Definition~\ref{def:component}\ref{en:StrongConn}.

\section{Proof of Theorem~\ref{thm:main}\ref{en:MainSatisfiable}}\label{sec:ProofSatisfiable}
The bulk of the proof will deal with the following setting.
\begin{prop}\label{prop:satisfiableSingleComponent}
Suppose that $U$ is an $L^3$-digraphon on $\OurSpace$. By Lemma~\ref{lem:LpImpliesBoundedCompact}, we have that $\rho_3(U)$ is defined. If $\rho_3(U)<1$, then for the random variable $N_n$ counting the number of bicycles in $\G(n,U)$, we have $\Expectation[N_n]\to 0$ as $n\to\infty$.
\end{prop}
\subsection{Proposition~\ref{prop:satisfiableSingleComponent} implies Theorem~\ref{thm:main}\ref{en:MainSatisfiable}} 
Let $U:=\overrightarrow{W}\left\llbracket \bigcup_{i\in I^*}\Omega_i\right\rrbracket$. We have $U\in L^3(\OurSpace^2)$ and the associated integral kernel operator $U$ is bounded and compact as an operator $U:L^{3}(\OurSpace)\to L^{3}(\OurSpace)$. Restricting our attention to $\bigcup_{i\in I^*}\Omega_i$ is sensible as contradictory cycles cannot appear elsewhere. This is formally stated in the next lemma.
\begin{lem}\label{lem:insideOne}
Suppose that $\phi\sim \TwoSAT(n,W)$ where in the sampling, we generated
$\tau:\Lit_n \to \Omega$. Then almost surely for each contradictory cycle
$u_1,\ldots,u_r$ in $\ImplDig(\phi)$ we have that
$\tau(u_1),\ldots,\tau(u_r)\in \bigcup_{i\in I^*}\Omega_i$.
\end{lem}
\begin{proof}
This is just Proposition~\ref{DIGRAPHONS-prop:cyclesconfined} in~\cite{HladkySavicky:Digraphons} translated to our setting.
\end{proof}
Let $M_n$ be the number of bicycles contained in some contradictory cycle in $\ImplDig(\TwoSAT(n,W))$. Let $N_n$ be the number of bicycles in $\G(n,U)$. Note that the edge set $F$ of each potential bicycle satisfies the conditions of Lemma~\ref{lem:marginalssame}. From Lemma~\ref{lem:marginalssame} and Lemma~\ref{lem:insideOne}, we get 
\begin{equation}\label{eq:ifthen}
\Probability[N_n=0]\le \Probability[M_n=0].
\end{equation}
We apply Proposition~\ref{prop:satisfiableSingleComponent} and get that $\Expectation[N_n]\to 0$, as $n\to \infty$. Markov's inequality and~\eqref{eq:ifthen} gives that $M_n=0$ asymptotically almost surely, as $n\to \infty$. By Proposition~\ref{prop:SATbicycle}, $\TwoSAT(n,W)$ is asymptotically almost surely satisfiable.

\subsection{Proof of Proposition~\ref{prop:satisfiableSingleComponent}}
We set up constants $\alpha$ and $C$ in a way which does not depend on $n$.
Let $\alpha>0$ be such that $\rho_3(U)+\alpha<1$. Next, we set up $C$. By Proposition~\ref{prop:GelfandLp}, we have $\rho_3(U)=\lim_{\ell\to\infty} \left(\|U^\ell\|_3\right)^{1/\ell}$. In particular, we can fix a constant $C>0$ such that for every $\ell\in \N_0$, we have
\begin{align}
\label{eq:CGelfand3}
\|U^\ell\|_3&\le C(\rho_3(U)+\alpha)^\ell\;.
\end{align}

For each $k \ge 2$, $2 \le a \le k$, and $1 \le b \le k-1$, let the random variable
$N^{k,a,b}_n$ count the number of $(k,a,b)$-bicycles in $\G(n,U)$. Let us first
focus on the path $u_1u_2\ldots u_k$ of any such bicycle. It is part of the
definition of a bicycle that the variables in the literals $u_i$ are pairwise distinct.
By symmetry, the probability that any such sequence forms a basis of a
$(k,a,b)$-bicycle is the same as for the sequence of positive literals
$v_1,\ldots,v_k$. That is, we have
\begin{align}
\nonumber
\Expectation\left[N^{k,a,b}_n\right]&=\binom{n}{k}k! \cdot 2^k \cdot
\Probability_{\G(n,U)}\left[\text{$v_1,\ldots,v_k$ forms a basis of a $(k,a,b)$-bicycle}\right]\\
\label{eq:EN}
&\le n^k\cdot 2^k\cdot\Probability_{\G(n,U)}\left[\text{$v_1,\ldots,v_k$ forms a basis of a $(k,a,b)$-bicycle}\right] \;.
\end{align}
In order for $v_1,\ldots,v_k$ to form a basis of a $(k,a,b)$-bicycle, the sequence
$v_1,\ldots,v_k$ has to form a path, and
additionally, two extra edges, namely $(\neg v_a, v_1)$ and $(v_k, \neg v_b)$
have to be present. 

We distinguish cases $a\ge b$ and $a<b$. 

\subsubsection*{Case $a\ge b$.}
We have
\begin{align}
\begin{split}\label{eq:MontPyth}
&\Probability_{\G(n,U)}\left[\text{$v_1,\ldots,v_k$ forms a basis of a $(k,a,b)$-bicycle}\right]
\\
&\le
\frac{1}{(2n)^{k+1}}\cdot \int_{x_a,x_b}U^b(\neg x_a,x_b)\cdot U^{a-b}(x_b,x_a) \cdot U^{k-a+1}(x_a,\neg x_b)
\;.
\end{split}
\end{align}
We use Hölder's inequality (HI) with exponents $\frac32$ and $3$, and the Cauchy--Schwarz Inequality (CSI),
\begin{align*}
&\int_{x_a,x_b}U^b(\neg x_a,x_b)\cdot U^{a-b}(x_b,x_a) \cdot U^{k-a+1}(x_a,\neg x_b)
\\
&\leBy{(HI)}
\left(\int_{x_a,x_b}\left(U^b(\neg x_a,x_b)\right)^{3/2}\cdot \left(U^{a-b}(x_b,x_a)\right)^{3/2}\right)^{2/3}
\cdot
\left(\int_{x_a,x_b}\left(U^{k-a+1}(x_a,\neg x_b)\right)^3\right)^{1/3}
\\
&\leBy{(CSI)}
\left(\int_{x_a,x_b}\left(U^b(\neg x_a,x_b)\right)^{3}\right)^{1/3}
\cdot
\left(\int_{x_a,x_b}\left(U^{a-b}(x_b,x_a)\right)^{3}\right)^{1/3}
\cdot
\left(\int_{x_a,x_b}\left(U^{k-a+1}(x_a,\neg x_b)\right)^3\right)^{1/3}
\\
&= \| U^b \|_3 \cdot \|U^{a-b}\|_3 \cdot \| U^{k-a+1} \|_3
\\
&\leByRef{eq:CGelfand3} C^3 (\rho_3(U)+\alpha)^{k+1}\;.
\end{align*}
We substitute this into~\eqref{eq:MontPyth},
\begin{equation*}
\Probability_{\G(n,U)}\left[\text{$v_1,\ldots,v_k$ forms a basis of a $(k,a,b)$-bicycle}\right]
\le
C^3\left(\frac{\rho_3(U)+\alpha}{2n}\right)^{k+1}\;,
\end{equation*}
which can in turn be substituted into~\eqref{eq:EN},
\begin{equation}\label{eq:UzS1}
\Expectation\left[N^{k,a,b}_n\right]\le \frac{C^3}{2n}\cdot \left(\rho_3(U)+\alpha\right)^{k+1}\;.
\end{equation}

\subsubsection*{Case $a< b$.}
The calculations are similar. We have
\begin{align*}
\begin{split}
&\Probability_{\G(n,U)}\left[\text{$v_1,\ldots,v_k$ forms a basis of a $(k,a,b)$-bicycle}\right]
\\
&\le
\frac{1}{(2n)^{k+1}}\cdot \int_{x_1,x_a,x_b,x_k}U(\neg x_a,x_1)U^{a-1}(x_1,x_a)\cdot U^{b-a}(x_a,x_b) \cdot U^{k-b}(x_b,x_k)U(x_k,\neg x_b)\;.
\end{split}
\end{align*}
We use Hölder's inequality (HI) with exponents $\frac32$ and $3$,
\begin{align*}
&\int_{x_1,x_a,x_b,x_k}U(\neg x_a,x_1)U^{a-1}(x_1,x_b) \cdot U^{k-b}(x_b,x_k)U(x_k,\neg x_b)\cdot U^{b-a}(x_a,x_b)
\\
&\leBy{(HI)}
\left(\int_{x_1,x_a,x_b,x_k}\left(U(\neg x_a,x_1)U^{a-1}(x_1,x_a)\right)^{3/2}\cdot \left(U^{k-b}(x_b,x_k)U(x_k,\neg x_b)\right)^{3/2}\right)^{2/3}
\cdot
\left(\int_{x_a,x_b}\left(U^{b-a}(x_a,x_b)\right)^3\right)^{1/3}
\\
&=
\left(\int_{x_1,x_a}\left(U(\neg x_a,x_1)U^{a-1}(x_1,x_a)\right)^{3/2}\right)^{2/3}
\cdot
\left(\int_{x_1,x_a}\left(U^{k-b}(x_b,x_k)U(x_k,\neg x_b)\right)^{3/2}\right)^{2/3}
\cdot
\|U^{b-a}\|_3
\;.
\end{align*}
We can now use the Cauchy--Schwarz inequality on the first term, 
\begin{align*}
\left(\int_{x_1,x_a}\left(U(\neg x_a,x_1)U^{a-1}(x_1,x_a)\right)^{3/2}\right)^{2/3}
&\le 
\left(\int_{x_1,x_a}\left(U(\neg x_a,x_1)\right)^{3}\right)^{1/3}
\cdot
\left(\int_{x_1,x_a}\left(U^{a-1}(x_1,x_a)\right)^{3}\right)^{1/3}\\
&=\| U \|_3 \cdot \| U^{a-1} \|_3
\;.    
\end{align*}
We can now use the Cauchy--Schwarz inequality on the second term in the same way. Combined, we conclude that
\begin{align*}
&\int_{x_1,x_a,x_b,x_k}U(\neg x_a,x_1)U^{a-1}(x_1,x_b) \cdot U^{k-b}(x_b,x_k)U(x_k,\neg x_b)\cdot U^{b-a}(x_a,x_b)\\
&\le \| U \|_3 \cdot \| U^{a-1} \|_3 \cdot \| U \|_3 \cdot \| U^{k-b } \|_3 \cdot \| U^{b-a} \|_3
\\
&\leByRef{eq:CGelfand3} C^5 (\rho_3(U)+\alpha)^{k+1}\;.
\end{align*}
This allows us to get a counterpart to~\eqref{eq:UzS1},
\begin{equation}\label{eq:UzS2}
\Expectation\left[N^{k,a,b}_n\right]\le \frac{C^5}{2n}\cdot \left(\rho_3(U)+\alpha\right)^{k+1}\;.
\end{equation}

\subsubsection*{Putting it together.}
We have
\begin{align*}
    \Expectation[N_n]\le \sum_{k=2}^\infty\sum_{a=2}^k\sum_{b=1}^{k-1} \Expectation\left[N^{k,a,b}_n\right]
    \;
    \leBy{by~\eqref{eq:UzS1},\eqref{eq:UzS2}}
 \frac{1}{2n} (C^3+C^5) \sum_{k=2}^\infty k^2 (\rho_3(U)+\alpha)^{k+1}\;.
\end{align*}
Since $\rho_3(U)+\alpha < 1$, we obtain $\Expectation[N_n] = O(\frac{1}{n})$.
This proves Proposition~\ref{prop:satisfiableSingleComponent}.

\section{Proof of Theorem~\ref{thm:main}\ref{en:MainUnsatisfiable}}\label{sec:Unsatisfiable}
First, we use Proposition~\ref{prop:approximatefrombelow} to find a graphon $V\le W$ with $\rho^*(V)>1$.

Let $\Omega_i$ be an arbitrary contradictory component with $\eta:=\rho\left(\overrightarrow{V}\llbracket \Omega_i\rrbracket\right)>1$. Let $\OurSpace^*:=\Omega_i$, and let $\OurMeasure^*$ be a measure on $\OurSpace^*$ naturally rescaled so that it becomes a probability measure, $\OurMeasure^*(A):=\frac{\OurMeasure(A)}{\OurMeasure(\Omega_i)}$ for every $A\subset \Omega_i$. By Proposition~\ref{prop:contradictorycomponentsproduct}, we have $\OurSpace^*=_0 \Lambda^*\times \signum$ and $\OurMeasure^*=\lambda^*\times \mu^{+-}$ for a certain subspace $\Lambda^*\subset \Lambda$ and similarly rescaled probability measure $\lambda^*$ on $\Lambda^*$. The negation map $\neg:\OurSpace \to \OurSpace$ from Section~\ref{ssec:OriginalSat} naturally restricts to $\neg:\OurSpace^* \to \OurSpace^*$.

We define versions of $V$ and $\overrightarrow{V}$ zoomed in to $\OurSpace^*$,
\begin{align}
\label{def:U}
&U\in L^\infty\left( (\OurSpace^*)^2\right)   \quad &U(x,y):=V(x,y)\;,
\\
\label{def:Gamma}
&\Gamma\in L^\infty\left( (\OurSpace^*)^2\right)   \quad &\Gamma(x,y):=\overrightarrow{V}(x,y)
\;.
\end{align}

In the remainder of this section, we will work (unless otherwise stated) within the space $(\OurSpace^*,\OurMeasure^*)$. In particular, note that $\rho(\Gamma)=\frac{\eta}{\OurMeasure(\Omega_i)}$.

Of course, we expect that there is a reason for unsatisfiability of $\OriginalTwoSAT\left(n,W\right)$ within $\OurSpace^*$. This is expressed in Proposition~\ref{prop:MainUnsat} below. Let us do some preparations to state it. For integers $a, b \ge 2$, an \emph{$(a,b)$-snake} is any formula of the form
\begin{equation} \label{eq:snake-F}
\begin{array}{ll}
F = & (\neg f \vee l_1) \wedge (\neg l_1 \vee l_2) \wedge \ldots \wedge (\neg l_{a-1} \vee \neg f) \wedge \\
    & (f \vee l_a) \wedge (\neg l_a \vee l_{a+1}) \wedge \ldots \wedge (\neg l_{a+b-2} \vee f)
\end{array}
\end{equation}
such that the variables of the literals $f,l_1, \ldots, l_{a+b-2}$ are all distinct. We emphasize that equality of two snakes can include reordering of the clauses, for example
\begin{align}
\begin{split}
\label{eq:4snakes}
&(\neg v_{99} \vee v_1) \wedge (\neg v_1 \vee v_2) \wedge (\neg v_2 \vee \neg v_{99}) \wedge (v_{99} \vee v_3) \wedge (\neg v_3 \vee v_{99})    \\
=&(\neg v_{99} \vee \neg v_2) \wedge (\neg(\neg v_2) \vee \neg v_1) \wedge (\neg(\neg v_1) \vee \neg v_{99}) \wedge (v_{99} \vee v_3) (\neg v_3 \vee v_{99})   \\
=&(\neg v_{99} \vee v_1) \wedge (\neg v_1 \vee v_2) \wedge (\neg v_2 \vee \neg v_{99}) \wedge (v_{99} \vee \neg v_3) \wedge (\neg(\neg v_3) \vee \neg v_{99})  \\
=&(\neg v_{99} \vee \neg v_2) \wedge (\neg(\neg v_2) \vee \neg v_1) \wedge (\neg(\neg v_1) \vee \neg v_{99}) \wedge (v_{99} \vee \neg v_3) \wedge (\neg(\neg v_3) \vee \neg v_{99})      
\end{split}
\end{align}
are~4 different choices of the literals $f,l_1, \ldots, l_{a+b-2}$ which lead to the same $(3,2)$-snake, and no more choices exist. Later, we will take numbers $a$ and $b$ growing with $n$ such that $a\neq b$ (see~\eqref{EQ:aNEQb}). Then in general, there are exactly~4 distinct choices of the literals $f,l_1, \ldots, l_{a+b-2}$ as in~\eqref{eq:snake-F} that lead to the same $(a,b)$-snake.

For an $(a,b)$-snake $F$, let \index{$\Lit(F)$}$\Lit(F)$ be the set of literals that have
an occurrence in $F$. Note that $|\Lit(F)| = 2(a+b)-2$, since $\Lit(F)$
is closed under negation, that is, if it contains a literal $\ell$, it
contains also $\neg \ell$.

When $a$ and $b$ are not specified, we call $F$ simply a \emph{snake}.
The definition of snakes appeared previously (with minor modifications) in literature concerning 2-SAT (e.g.~\cite{FriedrichRothenberger2022}). The importance of this definition is that any formula containing a snake as a subformula is obviously unsatisfiable.
\begin{prop}\label{prop:MainUnsat}
Let $U$ and $\eta$ be as above. Let $\beta\in (\frac{\OurMeasure(\Omega_i)}{\eta},\infty)$ be arbitrary. As $n\to\infty$, $\OriginalTwoSAT\left(\lfloor\beta n\rfloor,U\right)$ asymptotically almost surely contains a snake as a subformula.
\end{prop}

\subsection{Proposition~\ref{prop:MainUnsat} implies Theorem~\ref{thm:main}\ref{en:MainUnsatisfiable}}
Below, we shall argue that the model $\OriginalTwoSAT(n,V)$ is asymptotically almost surely unsatisfiable. Observe that this will imply the desired asymptotic almost sure unsatisfiability of $\OriginalTwoSAT(n,W)$. Indeed, as $V\le W$, the random clauses of $\OriginalTwoSAT(n,V)$ are stochastically dominated by the random clauses of $\OriginalTwoSAT(n,W)$.

Set $\beta:=\frac{2\OurMeasure(\Omega_i)}{1+\eta}$. The point of this choice is that we have $\beta\in (\frac{\OurMeasure(\Omega_i)}{\eta},\infty)$ (and so Proposition~\ref{prop:MainUnsat} applies) and $\beta<\OurMeasure(\Omega_i)$ (which will be important when we apply the law of large numbers below).
We take the asymptotics $n\to\infty$. 
Consider the way a random formula $\OriginalTwoSAT(n,V)$ was generated in Section~\ref{ssec:OriginalSat}. Let $I\subset [n]$ be the set of indices $i$ such that $x_i\in \Lambda^*$. For every $J\subset [n]$ of size $\lfloor \beta n\rfloor$, let $\mathcal{E}_J$ be the event that $J\subset I$, and that $|I\cap\{1,2,\ldots,\max(J)\}|=|J|$. That is, $\mathcal{E}_J$ is the event that at least $\lfloor \beta n\rfloor$ many elements $x_i$  were sampled from $\Lambda^*$, and that those $\lfloor \beta n\rfloor$ many one with the smallest indices $i$ form the set $J$. Define the event $\mathcal{O}$ defined by $|I|<\lfloor \beta n\rfloor$. We have that the events $\mathcal{O}$ and $\{\mathcal{E}_J\}_J$ partition our probability space. By the law of large numbers, asymptotically almost surely, $|I|=(1+o(1))\lambda(\Lambda^*)n=(1+o(1))\OurMeasure(\Omega_i)n$. In particular, $\Probability[\mathcal{O}]=o(1)$.

For any $J$ as above, the formula $\OriginalTwoSAT(n,V)$ in the conditional space $\mathcal{E}_J$ contains (up to renaming the variables) a formula $\OriginalTwoSAT\left(\lfloor\beta n\rfloor,U\right)$. By Proposition~\ref{prop:MainUnsat}, the latter formula contains a snake asymptotically almost surely, and thus is unsatisfiable asymptotically almost surely. We have
\begin{align*}
\Probability\left[\mbox{$\OriginalTwoSAT(n,W)$ is satisfiable}\right]&\le 
\Probability[\mathcal{O}]+\sum_{J\in \binom{[n]}{\lfloor \beta n\rfloor}}
\Probability[\mathcal{E}_J]\cdot\Probability\left[\mbox{$\OriginalTwoSAT(n,W)$ is satisfiable}\:|\: \mathcal{E}_J\right]\\
&\le \Probability[\mathcal{O}]+\sum_{J}
\Probability[\mathcal{E}_J]\cdot\Probability\left[\mbox{$\OriginalTwoSAT(\lfloor \beta n\rfloor,U)$ is satisfiable}\right]\\
&\le \Probability[\mathcal{O}]+\left(\sum_{J}
\Probability[\mathcal{E}_J]\right)\cdot\Probability\left[\mbox{$\OriginalTwoSAT(\lfloor \beta n\rfloor,U)$ is satisfiable}\right]
\le o(1)+1 \cdot o(1)\;.
\end{align*}

\subsection{Proof of Proposition~\ref{prop:MainUnsat}}
In the proof, we treat the number $N:=\lfloor \beta n\rfloor$ as $N=\beta n$. Let $\rho:=\frac{\eta}{\OurMeasure(\Omega_i)}=\rho(\Gamma)$.

Let $D$ be the peripheral multiplicity of $\Gamma$. Fix a partition $\OurSpace^*=X_0\sqcup X_1\sqcup \ldots\sqcup X_{D-1}$ as in Definition~\ref{def:graphicallyperiodic} provided by Theorem~\ref{thm:periodicity}. Fix two indices $\iota_{\mathrm{source}},\iota_{\mathrm{sink}}\in \{0,1,\ldots,D-1\}$ such that 
\begin{equation}\label{eq:contradictorysourcesink}
\lambda^*\left(\{y\in \Lambda^*:(y,+)\in X_{\iota_{\mathrm{source}}}\mbox{ and }(y,-)\in X_{\iota_{\mathrm{sink}}}\}\right)>0\;.    
\end{equation}

\subsubsection{Setup for counting snakes}
We first introduce an abstract version of our setup. That is, first we work in a general probability space. Suppose that $\{Z_s\}_{s\in \mathcal{S}}$ are indicator random variables indexed by a finite set $\mathcal{S}$. Let $Z=\sum_{s\in \mathcal{S}} Z_s$ be their sum. 
Let $\mathcal{T}\subset\mathcal{S}^2$ be an arbitrary superset of the set $\{(s_1,s_2)\in \mathcal{S}^2:s_1\neq s_2, \Expectation[Z_{s_1}Z_{s_2}]\neq \Expectation[Z_{s_1}]\Expectation[Z_{s_2}]\}$. 
The Paley--Zygmund inequality tells us that 
\begin{align}
\nonumber
\Probability[Z=0]\le \frac{\Expectation[Z^2]-\Expectation[Z]^2}{\Expectation[Z^2]}
&=
\frac{\sum_{(s_1,s_2)\in \mathcal{S}^2}(\Expectation[Z_{s_1}Z_{s_2}]-\Expectation[Z_{s_1}]\Expectation[Z_{s_2}])}{\Expectation[Z^2]}
\\
\nonumber
\JUSTIFY{we have $Z_s^2=Z_s$}&=
\frac{\sum_{s\in \mathcal{S}}(\Expectation[Z_{s}]-\Expectation[Z_{s}]^2)
+
\sum_{(s_1,s_2)\in \mathcal{T}}(\Expectation[Z_{s_1}Z_{s_2}]-\Expectation[Z_{s_1}]\Expectation[Z_{s_2}])}{\Expectation[Z^2]}
\\
\nonumber
&\le
\frac{\sum_{s\in \mathcal{S}}\Expectation[Z_{s}]
+
\sum_{(s_1,s_2)\in \mathcal{T}}\Expectation[Z_{s_1}Z_{s_2}]}{\Expectation[Z^2]}\\
\label{eq:second-moment:2}
&
=
\frac{\Expectation[Z]
+
\sum_{(s_1,s_2)\in \mathcal{T}}\Expectation[Z_{s_1}Z_{s_2}]}{\Expectation[Z^2]}
\;.
\end{align}

We shall use the above setup for counting $(a,b)$-snakes on the set of variables $\{v_1,\ldots,v_N\}$, where 
\begin{equation}\label{eq:defAB}
a:=2D\cdot \lceil \log_{\rho} N\rceil + \iota_{\mathrm{sink}}-\iota_{\mathrm{source}}    
\;\mbox{ and }\;
b:=3D\cdot \lceil \log_{\rho} N\rceil + \iota_{\mathrm{source}}-\iota_{\mathrm{sink}}
\;.
\end{equation}
Note that 
\begin{equation}\label{EQ:aNEQb}
    a\neq b\;
\end{equation}
for $N$ large enough. 

That is, let $\mathcal{S}$ be the set of snakes on variables $\{v_1,\ldots,v_N\}$. Taking into account the discussion around~\eqref{eq:4snakes}, we have $|\mathcal{S}|=2^{a+b-3}(N)_{a+b-1}$,
where $(N)_{a+b-1}$ is a falling factorial. If $a+b=o(\sqrt{N})$, then this gives 
\begin{equation}\label{eq:num-snakes}
|\mathcal{S}| = \Theta((2N)^{a+b-1})\;.
\end{equation}

For every $F \in \mathcal{S}$, let $Z_F$ be an indicator variable
representing the occurrence of $F$ in our random formula $\OriginalTwoSAT (N,U )$ obtained in the equivalent model $\TwoSAT (N,U )$ as described in Section~\ref{ssec:convenientmodel}.
Let $Z = \sum_{F \in \mathcal{S}} Z_F$. For a formula $F$, write $\Var(F)$ for the set of all variables it contains. We define
\begin{equation} \label{eq:define-D}
\mathcal{T} := \{(F, G)\in \mathcal{S}^2 \::\: F \not= G \mbox{ and } \Var(F) \cap \Var(G) \not= \emptyset\}\;.
\end{equation}
Obviously, this $\mathcal{T}$ satisfies the assumptions above.

An $(a,b)$-snake $F$ can be represented as a cycle of implications
in several ways.
Namely, $F$ is equivalent to each of the following $4$ different cycles
\begin{equation} \label{eq:four-cycles}
\begin{array}{l}
f \to l_1 \to \ldots \to l_{a-1} \to \neg f \to
l_a \to \ldots \to l_{a+b-2} \to f \\
f \to \neg l_{a-1} \to \ldots \to \neg l_1 \to \neg f \to
l_a \to \ldots \to l_{a+b-2} \to f \\
f \to l_1 \to \ldots \to l_{a-1} \to \neg f \to
\neg l_{a+b-2} \to \ldots \to \neg l_a \to f \\
f \to \neg l_{a-1} \to \ldots \to \neg l_1 \to \neg f \to
\neg l_{a+b-2} \to \ldots \to \neg l_a \to f
\end{array}
\end{equation}
If $C$ is a rooted cycle in any of the forms~(\ref{eq:four-cycles}) we call it an \emph{$(a,b)$-serpent}. We let
\index{$\mathrm{clauses}(C)$}
$\mathrm{clauses}(C)$ denote the snake $F$. Moreover, let $\impl(F)$
be the set of all implications equivalent to some of the clauses of $F$.
Hence, $\impl(F)$ is the union of the four serpents~(\ref{eq:four-cycles})
(this union is generated by the first and the fourth serpent of~\eqref{eq:four-cycles}, and also by the second and the third serpent of~\eqref{eq:four-cycles}).

\subsubsection{First and second moment for snakes}
The following lemma handles the probability of containment of a particular snake. Recall that $U$ is bounded (see~\eqref{def:U}).
\begin{lem} \label{lem:expected-F}
The probability that any given $(a,b)$-snake $F$ is contained in $\OriginalTwoSAT (N,U)$ is
\begin{equation}\label{eq:integralexpression}
\Probability[Z_F] = \frac{1}{(2N)^{a+b}}
\int_{x \in \OurSpace^*} \Gamma^a(x, \neg x) \Gamma^b(\neg x, x) \diff \OurMeasure^*(x) \;, \end{equation}
when $N \ge 2\| U\|_\infty$.
\end{lem}
\begin{proof}
The quantity $\Probability[Z_F]$ can be expressed using the way we generate random formulas in Section~\ref{ssec:convenientmodel}. That is, we integrate over the representatives $x, y_1, \ldots, y_{a+b-2}$ of the literals of the snake, and use~\eqref{eq:2SAT-using-tau} on every clause of the snake (see~\eqref{eq:snake-F}). This results in a $(a+b)$-fold product of terms coming from~\eqref{eq:2SAT-using-tau} in the integral,
\begin{align*}
\Probability[Z_F] 
= \frac{1}{(2N)^{a+b}} \int_{x, y_1, \ldots, y_{a+b-2} \in \OurSpace^*}& 
U(\neg x, y_1) U(\neg y_1, y_2) \ldots U(\neg y_{a-1}, \neg x) \cdot\\
&\cdot U(x, y_a) U(\neg y_a, y_{a+1}) \ldots U(\neg y_{a+b-2}, x) 
\diff \OurMeasure^*(x)\, \diff \OurMeasure^*(y_1) \ldots \diff \OurMeasure^*(y_{a+b-2})\;.
\end{align*}
Using the definition of $\Gamma$, this can be rewritten as
\begin{align*}
\Probability[Z_F] = \frac{1}{(2N)^{a+b}} \int_{x, y_1, \ldots, y_{a+b-2} \in \OurSpace^*} 
&\Gamma(x, y_1) \Gamma(y_1, y_2) \ldots \Gamma(y_{a-1}, \neg x) \cdot\\
&\cdot \Gamma(\neg x, y_a) \Gamma(y_a, y_{a+1}) \ldots \Gamma(y_{a+b-2}, x) 
\diff \OurMeasure^*(x)\, \diff \OurMeasure^*(y_1) \ldots \diff \OurMeasure^*(y_{a+b-2})\;.
\end{align*}
Note that each of the variables $y_1, \ldots, y_{a-1}$ has an occurrence
in exactly two terms in the integrand and these terms have the form
$\Gamma(z_1, y_i) \Gamma(y_i, z_2)$ with some other variables $z_1$
and $z_2$. Since no other term contains $y_i$, the integration over $y_i$
produces $\Gamma^2(z_1, z_2)$. Using this by induction, we can
integrate over all the variables $y_1, \ldots, y_{a-1}$ and replace the
product of the terms containing them by $\Gamma^a(x, \neg x)$.
Similarly, the terms containing $y_a, \ldots, y_{k-2}$ can be replaced by
$\Gamma^b(\neg x, x)$.
\end{proof}

Observe that Proposition~\ref{prop:asymptoticpowers} and Lemma~\ref{lem:eigenvectorsbounded} tell us that there exists a constant $c\ge 0$ such that for every $\ell\in \N$, 
\begin{equation} \label{eq:define-C}
\|\Gamma^\ell\|_\infty \le c \rho^\ell
\;.
\end{equation}

The next lemma uses Lemma~\ref{lem:expected-F} to get the order of magnitude of $\Probability[Z_F]$.
\begin{lem}\label{lem:orderFsnake}
For each $(a,b)$-snake $F$ we have
\begin{equation} \label{eq:expected-F-assumption}
\Probability[Z_F] = \Theta\left(\left( \frac{\rho}{2N} \right)^{a+b}\right) \;.
\end{equation}
\end{lem}
\begin{proof}
For the upper bound, we use~\eqref{eq:define-C} and see that the integrand in~\eqref{eq:integralexpression} is at most $c^2\rho^{a+b}$. Thus, $\Probability[Z_F]\le c^2\left( \frac{\rho}{2N} \right)^{a+b}$.

For the lower bound, we again invoke~\eqref{eq:integralexpression}, but this time together with a more careful estimate directly from Proposition~\ref{prop:asymptoticpowers}. Let $\alpha\in (0,\rho)$ be given by Proposition~\ref{prop:asymptoticpowers} for the digraphon $\Gamma$. Let $v_L$ and $v_R$ be the left and the right eigenvector for the eigenvalue $\rho(\Gamma)$ of $\Gamma$. By Theorem~\ref{thm:Schaefer74}, the functions $v_L,v_R\in L^2(\OurSpace^*)$ are positive almost everywhere. Combining with~\eqref{eq:contradictorysourcesink}, we see that there exists $\delta>0$ such that for the set
\[
R:=\{y\in \Lambda^*:(y,+)\in X_{\iota_{\mathrm{source}}}, v_L\big((y,+)\big)\ge \delta\mbox{ and }(y,-)\in X_{\iota_{\mathrm{sink}}}, v_R\big((y,-)\big)\ge \delta\}\;,
\]
we have $\lambda^*(R)\ge \delta$. To obtain a lower bound on the integral in~\eqref{eq:integralexpression}, we proceed as follows. We have
\begin{align}\label{eq:DSM}
\int_{x \in \OurSpace^*} \Gamma^a(x, \neg x) \Gamma^b(\neg x, x) \diff \OurMeasure^*(x)
\ge 
\frac{1}2\int_{y \in R} \Gamma^a\big((y,+), (y,-)\big) \Gamma^b\big((y,-), (y,+)\big) \diff \lambda^*(y)
\;.
\end{align}
For $y \in R$, we have $(y,+)\in X_{\iota_{\mathrm{source}}}$ and $(y,-)\in X_{\iota_{\mathrm{sink}}}$. By~\eqref{eq:defAB}, we have that $a\equiv \iota_{\mathrm{sink}}-\iota_{\mathrm{source}} \mod D$ and $b\equiv \iota_{\mathrm{source}}-\iota_{\mathrm{sink}} \mod D$. Thus, we can continue with~\eqref{eq:DSM} with the help of the point estimate from Proposition~\ref{prop:asymptoticpowers},
\begin{align*}
&\int_{x \in \OurSpace^*} \Gamma^a(x, \neg x) \cdot \Gamma^b(\neg x, x) \diff \OurMeasure^*(x)\\
&\ge 
\frac{1}2\int_{y \in R} \left(\rho^a v_R\big((y,+)\big) v_L\big((y,-)\big)+O(\alpha^a)\right) \cdot \left(\rho^b v_R\big((y,-)\big) v_L\big((y,+)\big)+O(\alpha^b)\right) \diff \lambda^*(y)\\
\JUSTIFY{def of $R$}&\ge 
\frac{\delta}{2}\left(\rho^a \delta^2+O(\alpha^a)\right) \cdot \left(\rho^b \delta^2+O(\alpha^b)\right)
=
\Theta\left(\rho^{a+b}\right)\;.
\end{align*}
It suffices to plug this bound into~\eqref{eq:integralexpression}.
\end{proof}

\begin{lem} \label{lem:first-moment}
For the sum $Z:=\sum_{F\in\mathcal{S}} Z_F$ we have $\Expectation[Z] = \Theta\left( \frac{\rho^{a+b}}{N} \right)$. In particular, by the choice of $a$ and $b$ in~\eqref{eq:defAB}, we have $\Expectation[Z]\to \infty$.
\end{lem}
\begin{proof}
We combine~\eqref{eq:num-snakes} and Lemma~\ref{lem:orderFsnake}. 
\end{proof}

To use the bound~(\ref{eq:second-moment:2}), we define
$\Delta := \sum_{(F,G) \in \mathcal{T}} \Expectation[Z_F Z_G]$. To prove Proposition~\ref{prop:MainUnsat}, we will show that $\Delta=o(\Expectation[Z^2])$.
By symmetry, we can fix an arbitrary $F_0\in\mathcal{S}$. Below, we write
\begin{equation*}
F_0 =  (\neg f \vee l_1) \wedge (\neg l_1 \vee l_2) \wedge \ldots \wedge (\neg l_{a-1} \vee \neg f) \wedge
     (f \vee l_a) \wedge (\neg l_a \vee l_{a+1}) \wedge \ldots \wedge (\neg l_{a+b-2} \vee f)
\;.
\end{equation*}
We then have
\begin{equation} \label{eq:reformulate-delta}
\Delta = |\mathcal{S}| \sum_{G\in\mathcal{S}: (F_0,G) \in \mathcal{T}} \Expectation[Z_{F_0} Z_G]
\;.
\end{equation}
Every term $Z_{F_0} Z_G$ in~\eqref{eq:reformulate-delta} is the indicator of the event the concatenation of the formulas $F_0$ and $G$ appears in $\TwoSAT (N,U )$. We denote $Z_{F_0\cup G}:=Z_{F_0}Z_G$.

In order to get an upper bound on $\Expectation[Z_{F_0 \cup G}]$, we use a generalization of the method used to prove Lemma~\ref{lem:expected-F}.
Consider an arbitrary mapping $\psi:\Var(F_0 \cup G) \to \OurSpace^*$.
Extend it into a mapping $\psi:\Lit(F_0 \cup G) \to \OurSpace^*$
by defining $\psi(\neg v) = \neg \psi(v)$.
Let $t(F_0 \cup G, \Gamma, \psi)$ be the product of the weights
of all clauses in $F_0 \cup G$, where a clause $(h \vee h')$
has weight 
\begin{equation}\label{eq:TY}
\frac{U(\psi(h), \psi(h'))}{2N}
=
\frac{\Gamma(\psi(\neg h), \psi(h'))}{2N}
=
\frac{\Gamma(\psi(\neg h'),\psi(h))}{2N}\;.
\end{equation}
Moreover, we consider $\psi(v)$ for every $v \in \Var(F_0 \cup G)$ as
a variable whose range is $\OurSpace^*$ and $t(F_0 \cup G, \Gamma, \psi)$
as a function of these variables. 
The expectation $\Expectation[Z_{F_0 \cup G}]$ will be expressed
as an integral of $t(F_0 \cup G, \Gamma, \psi)$ over the variables
in $V$, so we write
\begin{equation} \label{eq:expected-FG-1}
\Expectation[Z_{F_0 \cup G}] = \int_{\psi} t(F_0 \cup G, \Gamma, \psi) \, \diff((\OurMeasure^*)^{\Var(F_0\cup G)}) \;.
\end{equation}

Given any serpent $C=(g \to l'_1 \to \ldots \to l'_{a-1} \to \neg g \to
l'_a \to \ldots \to l'_{a+b-2} \to g)$, we say that $C$ is \emph{original} if 
we have $\clauses(C)\neq F_0$. We say that $C$ is \emph{overlapping} if $(F_0,\clauses(C))\in\mathcal{T}$. Note that by~\eqref{eq:4snakes}, serpents are in a 4-to-1 correspondence to all $(a,b)$-snakes. That is, we can rewrite~\eqref{eq:reformulate-delta} as
\begin{equation} \label{eq:reformulate-deltausingserpents}
\Delta = \tfrac14|\mathcal{S}| \sum_{\textrm{$C$ original overlapping serpent}} \Expectation[Z_{F_0\cup C}]
\;.
\end{equation}
(Here, $Z_{F_0\cup C}$ is a short for $Z_{F_0\cup \clauses(C)}$. Similarly, later on, $\Var(C)$ will be a short for $\Var(\clauses(C))$.) 

\begin{figure}\centering
	\includegraphics[scale=0.9]{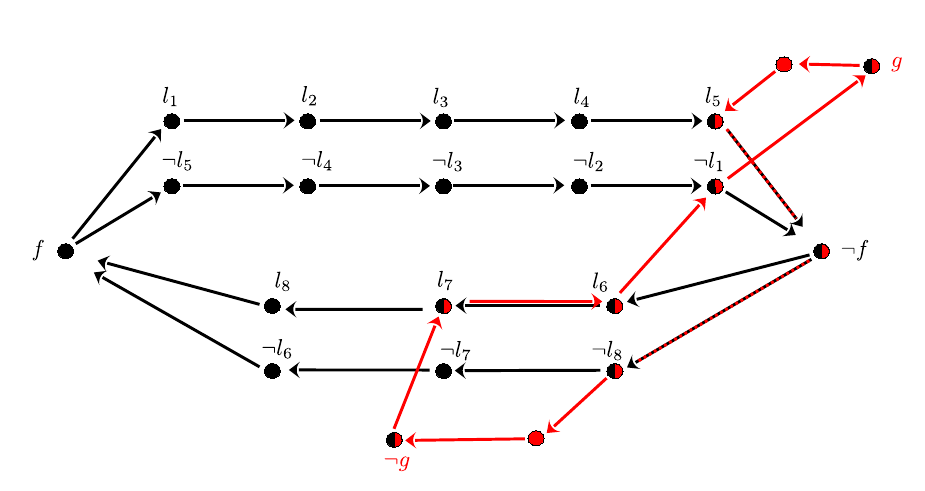}
	\caption{An example of a serpent when $a=6$, $b=4$. The serpent is depicted in red. Its free literals are depicted in full red circles, its non-free literals are depicted in half-red circles. Its intersection edges are dotted. The good non-intersection sequences written in the format $(start,end,length)$ are: $(g,l_5,2)$, $(\neg l_8,\neg g,2)$, $(\neg g,l_7,1)$.}
	\label{fig:serpent}
\end{figure}

To get an upper bound on $\Delta$, we classify serpents, and obtain bounds on $\Expectation[Z_{F_0\cup C}]$ depending on the particular classification of serpent $C$. This is done in the following definition. See also Figure~\ref{fig:serpent} for an illustration.

\begin{defi}
Suppose that $C$ is an original serpent rooted at a literal $g$. A literal in $C$ is called \emph{free}, if it is not in $\Lit(F_0) \cup \{g, \neg g\}$. A directed edge of $C$ is an \emph{intersection edge} if it belongs to $C \cap \impl(F_0)$, and a \emph{non-intersection edge} otherwise.
A nonempty sequence of consecutive (i.e., increasing along $C$) non-intersection edges of $C$, such
that the terminal vertices of the sequence are non-free literals and
there is no other non-free literal in the sequence, is called a \emph{good non-intersection sequence}. The \emph{length} of a good non-intersection sequences is~1 plus the number of free literals in it.

\index{$r(\mathcal{C})$,$s(\mathcal{C})$}Let $r(C)$ be the number of free literals in $C$. Let $s(C) = |C \setminus \impl(F_0)|$.

We say that $C$ is of \emph{Type~I} if $g \not\in \Lit(F_0)$. We say that $C$ is of \emph{Type~II} if $g \in \Lit(F_0) \setminus \{f, \neg f\}$. We say that $C$ is of \emph{Type~III} if $g \in \{f, \neg f\}$.
\end{defi}

\begin{defi}\label{def:EquivalentSerpents}
Suppose that $C$ is a serpent. We define the \emph{intersection pattern} of $C$ as the rooted cycle coming from $C$ (with the same root) where we replace the free literals by positive literals on new variables $z_1, z_2, \ldots, z_{r(C)}$ (not appearing in $\Lit_N$) in the order defined by the numbering of the literals in $C$. Hence, the intersection pattern is uniquely determined by $C$. 

For a given serpent $C$, let \index{$[C]$}$[C]$ be the equivalence class of all serpents having the same intersection pattern, and let \index{$\mathfrak{C}$}$\mathfrak{C}$ be the set of all these equivalence classes.

Note that all serpents $C$ within one intersection pattern have the same type (Type~I/Type~II/Type~III) and the same parameters $r(C)$ and $s(C)$. Thus we can define \index{$r(\mathcal{E})$,$s(\mathcal{E})$}$r(\mathcal{E})$ and $s(\mathcal{E})$ for $\mathcal{E}\in\mathfrak{C}$.
\end{defi}

\begin{lem} \label{lem:good-subseq}
Suppose that $C$ is an original serpent. The number of good non-intersection sequences of $C$ is $s(C)-r(C) \ge 1$. The sum of the lengths of all good non-intersection sequences of $C$ is equal to $s(C)$
\end{lem}
 
\begin{proof}
Let $g$ be the root of $C$. Since $g$ and $\neg g$ are non-free literals, every maximal sequence
of consecutive non-intersection edges and free literals has non-free
terminal vertices, so it is a good non-intersection sequence. The number
of edges in a good non-intersection sequence is~1 plus the number
of free literals in it. Hence, it contributes~1 to $s(C)-r(C)$. Moreover,
there is at least one non-intersection edge, since $F_0 \not= \clauses(C)$. This
implies that there is at least one good non-intersection sequence,
so we have $s(C)-r(C) \ge 1$.

The claim about the sum of the lengths is obvious.
\end{proof}

Note that a serpent has no free literals if it goes through the set
of literals $\Lit(F_0)$ in a way different from~$F_0$.

Let us derive an upper bound on $\Expectation[X_{F_0 \cup C}]$ in terms of the parameters $r(C)$ and $s(C)$.

\begin{lem}
Let $C$ be a serpent. We  have
\begin{equation} \label{eq:upper-bound-FG-2}
\Expectation[Z_{F_0 \cup C}] \le
\Theta\left(\left( \frac{\rho}{2N} \right)^{a+b}
c^{s(C)-r(C)}
\left(\frac{\rho}{2N}\right)^{s(C)}
\right)
\end{equation}
\end{lem}

\begin{proof}
Let $g$ be the root of $C$. Consider the following partition of the set of variables $\Var(F_0\cup C)=V_1\sqcup V_2\sqcup V_3$, where
\begin{align*}
V_1 &:= \Var(F_0) \;,\\
V_2 &:= \Var(\{g\}) \setminus \Var(F_0)\;,\\
V_3 &:= \Var(C) \setminus (\Var(\{g\}) \cup \Var(F_0))\;.
\end{align*}
Note that $V_2$ is either empty or a singleton and $V_3$
are the variables from $V$ corresponding to the free literals of $C$.
The integral (\ref{eq:expected-FG-1}) becomes
\begin{equation} \label{eq:expected-FG-2}
\Expectation[Z_{F_0 \cup C}] =
\int_{\psi_1\in (\OurSpace^*)^{V_1}} \int_{\psi_2\in (\OurSpace^*)^{V_2}} \int_{\psi_3\in (\OurSpace^*)^{V_3}} t(F_0 \cup C, \Gamma, \psi_1\oplus\psi_2\oplus \psi_3)
\; \diff((\OurMeasure^*)^{V_3}) \diff((\OurMeasure^*)^{V_2}) \diff((\OurMeasure^*)^{V_1})
\;.
\end{equation}

First, we prove that
\begin{equation} \label{eq:upper-bound-FG-1}
\Expectation[Z_{F_0 \cup C}] \le c^{s(C)-r(C)}
\left(\frac{\rho}{2N}\right)^{s(C)}
\cdot 
\int_{\psi_1\in (\OurSpace^*)^{V_1}}  t(F_0, \Gamma, \psi_1)
\,
\diff((\OurMeasure^*)^{V_1})\;.
\end{equation}
Indeed, the integrand in~(\ref{eq:upper-bound-FG-1})
is the product of the weights of the edges of $F_0$. This is
the part of $t(F_0 \cup C, \Gamma, \psi_1\oplus\psi_2\oplus \psi_3)$ that is not affected by integration over $V_3$ and $V_2$. We now need to incorporate the contribution of the non-intersection edges to $t(F_0 \cup C, \Gamma, \psi_1\oplus\psi_2\oplus \psi_3)$. We do it by looking at the good non-intersection sequences one-by-one. Consider each good non-intersection sequence, say consisting of $\ell$ edges, and use (\ref{eq:define-C}) and~\eqref{eq:TY}. After integrating over the $\ell-1$ free variables in its free literals, we see from this that such a sequence contributes at most $c\cdot \left(\frac{\rho}{2N}\right)^{\ell}$.

Finally, we multiply the above upper bounds over all good non-intersection sequences.
The number of these sequences is $s(C)-r(C)$ and the sum of their lengths is $s(C)$, so the product over all good sequences is the term
$$
c^{s(C)-r(C)} \left(\frac{\rho}{2N}\right)^{s(C)}
$$
in~(\ref{eq:upper-bound-FG-1}).

Let us now move from (\ref{eq:upper-bound-FG-1}) to~\eqref{eq:upper-bound-FG-2}. To this end, it suffices to combine~(\ref{eq:upper-bound-FG-1}) with~\eqref{eq:expected-F-assumption}.
\end{proof}

\begin{prop}
Let $C$ be an arbitrary serpent. Then
\begin{equation} \label{eq:upper-bound-FG-3}
|\mathcal{S}| \sum_{C' \in [C]} \Expectation[Z_{F_0 \cup C'}] \le
\Theta\left(\frac{1}{2N}
\left(\frac{c}{2N}\right)^{s(C)-r(C)} \rho^{a+b+s(C)}
\right)
\;.
\end{equation}
\end{prop}
\begin{proof}
Each serpent $C' \in [C]$ is given by a choice of a replacement of the free literals of $C$. We have at most $(2N)^{r(C)}$ choices of the free literals. We combine this with~\eqref{eq:num-snakes} and~\eqref{eq:upper-bound-FG-2} and get the desired bound.
\end{proof}

For every combination of $r$, $s$ which are admissible parameters $r(C)$ and $s(C)$ for serpents, let $M_{r,s,1}$
denote the number of intersection patterns of serpents of Type~I such that the parameters $r$, $s$,
have the given values. We define similarly $M_{r,s,2}$ and $M_{r,s,3}$ for serpents of Type~II and Type~III.

\begin{lem}
For any $r,s\in \N_0$, we have
\begin{align}
\label{eq:MType1}
M_{r,s,1} &\le 4N(2(a+b))^{3(s-r)} \;,\\
\label{eq:MType2}
M_{r,s,2} &\le 4(a+b)(2(a+b))^{3(s-r)} \;\mbox{, and}\\
\label{eq:MType3}
M_{r,s,3} &\le 4(2(a+b))^{3(s-r)}\;.
\end{align}
\end{lem}
\begin{proof}
First, we look at~\eqref{eq:MType1}. A specific intersection pattern of Type~I and with parameters $r$, $s$ is determined by the following choices
\begin{itemize}
\item
one of at most $2N$ possible literals for the root $g$,
\item
at most one of two choices `up'/`down' explained below,
\item
for each good non-intersection sequence, by a choice of
its two terminal vertices and length. The terminal vertices are elements of
$\Lit(F_0) \cup \{g, \neg g\}$ which has size $2(a+b)$. The length of
a good non-intersection sequence is an integer between $1$ and $a+b$. By Lemma~\ref{lem:good-subseq}, there are $s-r$ good non-intersection sequences. Hence, we have
at most $(2(a+b))^{2(s-r)} (a+b)^{s-r} \le (2(a+b))^{3(s-r)})$ possibilities.
\end{itemize}

\begin{figure}\centering
	\includegraphics[scale=0.6]{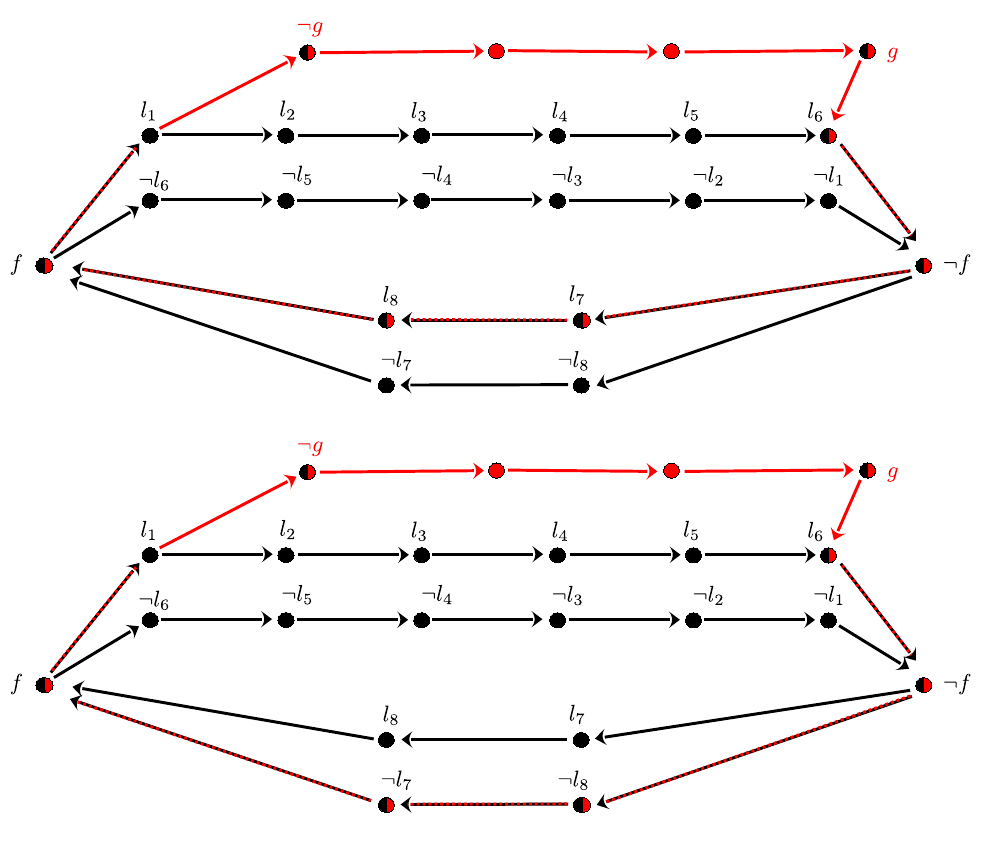}
	\caption{An example of two serpents which do not have the same intersection pattern but have the same set of good non-intersection sequences. Written in the format $(start,end,length)$, these are: $(g,l_6,1)$, $(l_1,\neg g,1)$, $(\neg g,g,3)$. Here, $a=7$, $b=3$.}
	\label{fig:serpentedges}
\end{figure}
Let us verify that the set of intersection edges is uniquely
determined by the above choices. These edges can be split
into maximum sequences of consecutive edges in $\impl(F_0)$. If $g$
is contained as an interior point of some of these sequences,
we split it further, so that $g$ is an terminal vertex. The number of
sequences of intersection edges is not important, so we can
safely do this. After this, the terminal vertices of intersection
sequences are determined by the root $g$
and the terminal vertices of good non-intersection sequences (that is, by the end of one good non-intersection sequence and the beginning of the next one). There is one exception, an example of which is shown in Figure~\ref{fig:serpentedges}. Namely, Figure~\ref{fig:serpentedges} shows two serpents, with two different intersection sequences, both starting at $l_6$ and ending at $l_1$. Such a situation can occur only for intersection sequences which entirely contain either one of the two branches of $F_0$ from $f$ to $\neg f$ or one of the two branches of $F_0$ from $\neg f$ to $f$. Each serpent contains at most one such intersection sequence. Hence, this non-uniqueness can be encoded by two choices, say `up'/`down'.

We now turn our attention to~\eqref{eq:MType2}. The argument is as for intersection patterns of Type~I, except that we can sharpen the bound on the number of choices of the root $g$ from $2N$ to $2(a+b)-2$.

Last, we cover~\eqref{eq:MType3}. The argument is again the same, with even a better bound on the number of choices of the root, namely $g \in \{f, \neg f\}$, i.e.,~2 choices.
\end{proof}

As we said, we prove Proposition~\ref{prop:MainUnsat} by using~\eqref{eq:second-moment:2}. For this, it remains to show that $\Delta=o(\Expectation[Z^2])$. To this end, we use~\eqref{eq:reformulate-deltausingserpents}. We partition $\mathfrak{C}=\mathfrak{C}_1\sqcup\mathfrak{C}_2\sqcup \mathfrak{C}_3$ into intersection patterns whose serpents are of Type~I, Type~II, or Type~III, respectively. That is, Proposition~\ref{prop:MainUnsat} will follow from the following three statements.
\begin{lem}\label{lem:TypeIaggregate}
We have
\[
|\mathcal{S}| \sum_{\mathcal{E}\in \mathfrak{C}_1}\sum_{C\in \mathcal{E}} \Expectation[Z_{F_0\cup C}]
=o(\Expectation[Z^2])\;.
\]
\end{lem}
\begin{lem}\label{lem:TypeIIaggregate}
We have
\[
|\mathcal{S}| \sum_{\mathcal{E}\in \mathfrak{C}_2}\sum_{C\in \mathcal{E}} \Expectation[Z_{F_0\cup C}]
=o(\Expectation[Z^2])\;.
\]
\end{lem}
\begin{lem}\label{lem:TypeIIIaggregate}
We have
\[
|\mathcal{S}| \sum_{\mathcal{E}\in \mathfrak{C}_3}\sum_{C\in \mathcal{E}} \Expectation[Z_{F_0\cup C}]
=o(\Expectation[Z^2])\;.
\]
\end{lem}

For the proofs of Lemmas~\ref{lem:TypeIaggregate}, \ref{lem:TypeIIaggregate}, \ref{lem:TypeIIIaggregate}, we partition $\mathfrak{C}_t$ (where $t=1,2,3$) as $\mathfrak{C}_t=\sqcup_{r,s} \mathfrak{C}_{s,s,t}$, where $r$ and $s$ are the parameters of the intersection pattern as in Definition~\ref{def:EquivalentSerpents}.

In fact, in the proofs of the three lemmas above, we shall prove that the left-hand sides are $o(\Expectation[Z]^2)$, which is a slightly stronger bound. To this end, we shall use Lemma~\ref{lem:first-moment}, which tells us that
\begin{equation}\label{eq:Tesco}
    \Expectation[Z]^2=\Theta\left(\frac{\rho^{2(a+b)}}{N^2}\right).
\end{equation}

\subsubsection{Proof of Lemma~\ref{lem:TypeIaggregate}}
Let $r\in \N_0$, $s\in \N$. Multiplying (\ref{eq:upper-bound-FG-3}) by~\eqref{eq:MType1}
yields
\begin{equation} \label{eq:upper-bound-t1-1}
|\mathcal{S}| \sum_{\mathcal{E}\in \mathfrak{C}_{r,s,1}}\sum_{C\in \mathcal{E}} \Expectation[Z_{F_0\cup C}]=
\Theta\left(
(2(a+b))^{3(s-r)}
\left(\frac{c}{2N}\right)^{s-r} \rho^{a+b+s}
\right) \;.
\end{equation}
We claim that $\mathfrak{C}_{r,s,1}=\emptyset$ for $s-r < 3$.
By Lemma \ref{lem:good-subseq}, $s-r$ is the number
of good non-intersection sequences.
In order to prove that we have at least~3 such sequences,
note that $C_G$ contains at least~3 non-free literals:
$g$, $\neg g$, and at least one further literal, say $h$ from $F_0$ (indeed, otherwise
we would have $\Var(F_0) \cap \Var(C) = \emptyset$ and $C$ would not be overlapping). Each of the literals $g$ and
$\neg g$ is connected to two non-intersection edges in $C$.
They contribute to $s-r$, but the contribution depends on whether
they belong to different good sequences. Two of them may belong
to the same sequence, if
they are in the branch between $g$ and $\neg g$ containing only free
literals. However, one of the branches contains $h$, so this branch
contributes at least $2$ to $s-r$ and the other at least $1$.

Now, fix $s\in \N$. The sum of (\ref{eq:upper-bound-t1-1}) over decreasing $r \ge 1$
satisfying $s-r \ge 3$ is a finite geometric progression with the
common ratio $(2(a+b))^3/(2N) = o(1)$. Hence, 
\[
|\mathcal{S}| \sum_{r\in \N: s-r\ge 3}\sum_{\mathcal{E}\in \mathfrak{C}_{r,s,1}}\sum_{C\in \mathcal{E}} \Expectation[Z_{F_0\cup C}]
=
\Theta\left(
(2(a+b))^{9}
\left(\frac{c}{2N}\right)^3 \rho^{a+b+s}
\right) \;.
\]
Since $\rho > 1$ is a constant, taking the sum over $1 \le s \le a+b$ yields
\begin{equation} \label{eq:upper-bound-t1-2}
|\mathcal{S}| \sum_{\mathcal{E}\in \mathfrak{C}_1}\sum_{C\in \mathcal{E}} \Expectation[Z_{F_0\cup C}]
=
\Theta\left(
(2(a+b))^{9}
\left(\frac{c}{2N}\right)^3 \rho^{2(a+b)}
\right) \;.
\end{equation}
To finish the proof of the lemma, we use~\eqref{eq:Tesco}. That is, we need to show the right-hand side of~\eqref{eq:upper-bound-t1-2} is $o\left(\frac{\rho^{2(a+b)}}{N^2}\right)$. By canceling  $\rho^{2(a+b)}$ on both sides, we see that this is equivalent to showing that 
\[
(2(a+b))^{9}
\left(\frac{c}{2N}\right)^3=o\left(\frac1{N^2}\right)
\;,
\]
which is evident.

\subsubsection{Proof of Lemma~\ref{lem:TypeIIaggregate}}
Let $r\in \N_0$, $s\in \N$. Multiplying (\ref{eq:upper-bound-FG-3}) by~\eqref{eq:MType2}
yields
\begin{equation} \label{eq:upper-bound-t2-1}
|\mathcal{S}| \sum_{\mathcal{E}\in \mathfrak{C}_{r,s,2}}\sum_{C\in \mathcal{E}} \Expectation[Z_{F_0\cup C}]
=
\Theta\left(
\frac{a+b}{N}
(2(a+b))^{3(s-r)}
\left(\frac{c}{2N}\right)^{s-r} \rho^{a+b+s}
\right) \;.
\end{equation}
We claim that $\mathfrak{C}_{r,s,2}=\emptyset$ for $s-r < 2$. By Lemma~\ref{lem:good-subseq}, we only need to treat the case $s-r = 1$. So, suppose $s-r = 1$. Then one of the halfcycles
of $C$, that is either the halfcycle starting at $g$ and ending
at $\neg g$ or the halfcycle starting at $\neg g$ and ending at $g$
consists only of intersection edges and due to the structure of
$\impl(F_0)$ contains both literals $f$ and $\neg f$. This is
a contradiction, since $C$ contains only one pair of contradictory
literals.

Now, fix $s\in \N$. Similarly as in the previous case, taking the sum of
(\ref{eq:upper-bound-t2-1}) over decreasing $r \ge 1$ satisfying
$s-r \ge 2$ yields
$$
|\mathcal{S}| \sum_{r\in \N_0: s-r\ge 2}\sum_{\mathcal{E}\in \mathfrak{C}_{r,s,2}}\sum_{C\in \mathcal{E}} \Expectation[Z_{F_0\cup C}]
=
\Theta\left(
\frac{a+b}{N}
(2(a+b))^{6}
\left(\frac{c}{2N}\right)^2 \rho^{a+b+s}
\right)
$$
and the sum over $1 \le s \le a+b$ is
\begin{equation} \label{eq:upper-bound-t2-2}
|\mathcal{S}| \sum_{\mathcal{E}\in \mathfrak{C}_2}\sum_{C\in \mathcal{E}} \Expectation[Z_{F_0\cup C}]
=
\Theta\left(
\frac{a+b}{N}
(2(a+b))^{6}
\left(\frac{c}{2N}\right)^2 \rho^{2(a+b)}
\right) \;.
\end{equation}
To finish the proof of the lemma, we use~\eqref{eq:Tesco}. That is, we need to show the right-hand side of~\eqref{eq:upper-bound-t2-2} is $o\left(\frac{\rho^{2(a+b)}}{N^2}\right)$. By canceling  $\rho^{2(a+b)}$ on both sides, we see that this is equivalent to showing that 
\[
\frac{a+b}{N}
(2(a+b))^{6}
\left(\frac{c}{2N}\right)^2
=o\left(\frac{1}{N^2}\right)
\;,
\]
which is evident.

\subsubsection{Proof of Lemma~\ref{lem:TypeIIaggregate}}
Let $r\in \N_0$, $s\in \N$.
Multiplying (\ref{eq:upper-bound-FG-3}) by~\eqref{eq:MType3}
yields
\begin{equation} \label{eq:upper-bound-t3-1}
|\mathcal{S}| \sum_{\mathcal{E}\in \mathfrak{C}_{r,s,3}}\sum_{C\in \mathcal{E}} \Expectation[Z_{F_0\cup C}]
=
\Theta\left(
\frac{1}{N}
(2(a+b))^{3(s-r)}
\left(\frac{c}{2N}\right)^{s-r} \rho^{a+b+s}
\right) \;.
\end{equation}
In order to get an upper bound, we distinguish two subcases,
namely $s-r \ge 2$ and $s-r=1$. We use the upper bound $M_{r,s,3}$
for each of these subcases, so we lose a factor of at most $2$.

Now, fix $s\in \N$. Similarly as in the previous case, taking the sum of
(\ref{eq:upper-bound-t3-1}) yields
$$
|\mathcal{S}| \sum_{r\in \N: s-r\ge 2}\sum_{\mathcal{E}\in \mathfrak{C}_{r,s,3}}\sum_{C\in \mathcal{E}} \Expectation[Z_{F_0\cup C}]
=
\Theta\left(
\frac{1}{N}
(2(a+b))^{6}
\left(\frac{c}{2N}\right)^2 \rho^{a+b+s}
\right) \;.
$$
and the sum over $1 \le s \le a+b$ is
\begin{equation*}
|\mathcal{S}| \sum_{s\in \N}\sum_{r\in \N_0: s-r\ge 2}\sum_{\mathcal{E}\in \mathfrak{C}_{r,s,3}}\sum_{C\in \mathcal{E}} \Expectation[Z_{F_0\cup C}]
=
\Theta\left(
\frac{1}{N}
(2(a+b))^{6}
\left(\frac{c}{2N}\right)^2 \rho^{2(a+b)}
\right) \;.
\end{equation*}
We immediately see that this term is negligible compared to~\eqref{eq:Tesco}. It remains to get the same conclusion about the terms with $s-r=1$.
If $s-r=1$, then~\eqref{eq:upper-bound-t3-1} gives
$$
|\mathcal{S}| \sum_{r\in \N: s-r=1}\sum_{\mathcal{E}\in \mathfrak{C}_{r,s,3}}\sum_{C\in \mathcal{E}} \Expectation[Z_{F_0\cup C}]
=
\Theta\left(
\frac{1}{N}
(2(a+b))^{3}
\left(\frac{c}{2N}\right) \rho^{a+b+s}
\right) \;.
$$
By Lemma~\ref{lem:good-subseq}, the above formula counts the contribution of intersection patterns with a single non-intersection sequence. One can verify that such
intersection patterns satisfy $s \le \max(a,b)=b$ (recall~\eqref{eq:defAB}).
Taking the sum over $1 \le s \le b$ yields
\begin{equation} \label{eq:upper-bound-t3-3}
|\mathcal{S}| \sum_{r\in \N: s-r=1}\sum_{\mathcal{E}\in \mathfrak{C}_{r,s,3}}\sum_{C\in \mathcal{E}} \Expectation[Z_{F_0\cup C}]
=
\Theta\left(
\frac{(\log N)^3}{N^2}
\rho^{a+2b}
\right)
\eqByRef{eq:Tesco}
O\left(\frac{(\log N)^3}{\rho^a}\cdot\Expectation[Z]^2\right)
\;,
\end{equation}
as was needed.

\section{Proof of Proposition~\ref{prop:rhoNula}}\label{sec:DenserRegimes}

\paragraph*{Equivalence $\ref{en:Dense1}\Leftrightarrow \ref{en:Dense2}$.} This is Proposition~\ref{DIGRAPHONS-prop:spectralradiusAndStrongComponents}\ref{DIGRAPHONS-en:CharSpectralRadius0} in~\cite{HladkySavicky:Digraphons}.

\paragraph*{Implication $\neg \ref{en:Dense1}\Rightarrow \neg \ref{en:Dense3}$.}
This is explained in Section~\ref{sssec:altarnitivescalings}, we just repeat the argument. Assume that $\rho^*(W)>0$. Let $f$ be an arbitrary function tending to infinity. We see that random clauses $\OriginalTwoSAT(n,f(n)\cdot W)$ eventually stochastically dominate random clauses of $\OriginalTwoSAT(n,U)$, where $U:=\frac{2W}{\rho^*(W)}$. Since we have $\rho^*(U)=\frac{2\rho^*(W)}{\rho^*(W)}>1$, Theorem~\ref{thm:main} tells us that $\OriginalTwoSAT(n,U)$ is asymptotically almost surely unsatisfiable. We conclude that $\OriginalTwoSAT(n,f(n)\cdot W)$ is asymptotically almost surely unsatisfiable, too.

\paragraph*{Implication $\ref{en:Dense4}\Rightarrow \ref{en:Dense3}$.} This is obvious.

\paragraph*{Implication $\ref{en:Dense2}\Rightarrow \ref{en:Dense4}$.} 
We will introduce another model of random 2-SAT, $\DensestTwoSAT(n,W)$. We will prove that assuming~\ref{en:Dense2}, 
\begin{enumerate}
    \item[(D1)] for every $n\in \N$ and for every $c>0$, random clauses of $\DensestTwoSAT(n,W)$ stochastically dominate random clauses of $\OriginalTwoSAT(n,c\cdot W)$, and 
    \item[(D2)] $\DensestTwoSAT(n,W)$ is asymptotically almost surely satisfiable.
\end{enumerate}
This will obviously prove that $\OriginalTwoSAT(n,f(n)\cdot W)$ is asymptotically almost surely satisfiable for every function $f$.
The model $\DensestTwoSAT(n,W)$ is defined exactly like $\OriginalTwoSAT(n,W)$ except that the insertion probability in~\eqref{eq:def2SATjb} is replaced by
\begin{equation*}
    \begin{cases}
    0& \mbox{if $W\left((x_i,\mathfrak{q}),(x_j,\mathfrak{s})\right)=0$,}\\
    1& \mbox{if $W\left((x_i,\mathfrak{q}),(x_j,\mathfrak{s})\right)>0$.}
    \end{cases}
\end{equation*}
It is obvious that~(D1) is satisfied. Let us now turn to~(D2). Let $\OurSpace=\Omega_0\sqcup \bigsqcup_{i\in I}\Omega_i$ be the decomposition of $\overrightarrow{W}$ into strong components. By  Proposition~\ref{DIGRAPHONS-prop:cyclesconfined} in~\cite{HladkySavicky:Digraphons}, we almost surely have the property that elements $y_{j_1},y_{j_2},\ldots,y_{j_k}\in \OurSpace$ representing an arbitrary cycle, say $C=v_{j_1}v_{j_2}\cdots v_{j_k}$, in the implication digraph $\ImplDig(\DensestTwoSAT(n,W))$ lie within one strong component, say $\Omega_i$. But as $\Omega_i$ is not a contradictory component, we have that $C$ is not a contradictory cycle, almost surely. We conclude that $\ImplDig(\DensestTwoSAT(n,W))$ contains no contradictory cycles almost surely. Thus, $\DensestTwoSAT(n,W)$ is satisfiable almost surely by Proposition~\ref{prop:SATcycle}.

\section{Further directions}
\subsection{Graphons with $\rho^*(W)=1$}\label{ssec:QuestionRhoOne}
If we take $W=\mathbbm{1}$, then the result of~\cite{MR1824274} implies that the probability of satisfiability of $\OriginalTwoSAT(n,\mathbbm{1})$ lies in the interval $(c,1-c)$ for some $c>0$. It is in fact reasonable to conjecture that the limit of $\Probability[\text{$\OriginalTwoSAT(n,\mathbbm{1})$ is satisfiable}]$ exists. More generally, it could be, that for every graphon $W$ with $\rho^*(W)=1$, we have that the limit of $\Probability[\text{$\OriginalTwoSAT(n,W)$ is satisfiable}]$ exists and lies in $(0,1)$. The main focus of~\cite{MR1824274} is on the `scaling window' of 2-SAT. This is a phenomenon first described in the setting of the giant component of random graphs, \cite{Boll:Window,Lucz:Window}. The main result of~\cite{MR1824274} asserts that 
\begin{align*}
0<\liminf_n \Probability[\mbox{$\OriginalTwoSAT(n,(1+Cn^{-1/3})\mathbbm{1})$ is satisf.}]\le \limsup_n \Probability[\mbox{$\OriginalTwoSAT(n,(1+Cn^{-1/3})\mathbbm{1})$ is satisf.}]<1
\end{align*}
for every $C\in \R$, and 
\begin{align*}
\liminf_n \Probability[\mbox{$\OriginalTwoSAT(n,(1-f(n))\mathbbm{1})$ is satisf.}]=1
\quad\mbox{and}\quad
\limsup_n \Probability[\mbox{$\OriginalTwoSAT(n,(1+f(n))\mathbbm{1})$ is satisf.}]=0\;,
\end{align*}
whenever $f:\N\to \R_+$ grows much faster than $n^{-1/3}$. The same question makes sense for any graphon with $\rho^*(W)=1$. Specifically, it is natural to ask whether the critical scaling is always of order $\Theta(n^{-1/3})$. This phenomenon --- often referred to as belonging to a `universality class' --- has been observed in many random discrete structures.

\subsection{Asymptotic rate of convergence in Theorem~\ref{thm:main}}
One might want to obtain the rate at which the probability of the satisfiability $\OriginalTwoSAT(n,W)$ approaches~0 or~1 in Theorem~\ref{thm:main}. The main result of~\cite{MR1824274} asserts that if $W\equiv C$ is a constant graphon, then the probability of satisfiability is exponentially small when $C>1$, and the probability of unsatisfiability is $\Theta(n^{-1})$ when $C\in (0,1)$. We believe that if $W$ is a graphon with $\rho^*(W)>1$ then $\Probability[\mbox{$\OriginalTwoSAT(n,W)$ is satisf.}]<\exp(-\eps_W n)$ for some $\eps_W>0$ that depends on $W$. The case $\rho^*(W)<1$ seems more delicate, allowing for different rates of the probability of satisfiability, depending on the structure of $W$.

\subsection{Higher clause sizes}
One could study inhomogeneous variants of random $k$-SAT for $k>2$. In the simplest version, the model would be given by a symmetric bounded measurable function $W:\OurSpace^k\to [0,\infty)$. Here, \emph{symmetric} means that $W(x_1,\ldots,x_k)=W(x_{\pi(1)},\ldots,x_{\pi(k)})$ for every permutation $\pi$. For $n\in \N$, we define random $k$-SAT
formula $\phi\sim\OriginalKSAT(n,W)$ on variables $\Lit_n$ by sampling elements
$x_1,\ldots,x_n\in \Lambda$ independently with distribution $\lambda$. For each
tuple $\{i_1,\ldots,i_k\}\in\binom{n}{k}$ and each tuple of signs $\mathfrak{s}_1,\ldots,\mathfrak{s}_k\in\signum$, we insert the clause
$\{\mathfrak{s}_1 v_{i_1}, \ldots,\mathfrak{s}_k v_{i_k}\}$ into $\phi$ independently from other choices with probability
\[
\min\left\{1,\frac{W\left((x_{i_1},\mathfrak{s}_1),\ldots,(x_{i_k},\mathfrak{s}_k)\right)}{n^{k-1}}\right\}\;.
\]
The scaling $\frac1{n^{k-1}}$ is chosen so that the resulting formula typically has $\Theta(n)$ clauses, which is the order of magnitude relevant for the homogeneous random $k$-SAT.

The theory of hypergraph limits (as established in~\cite{ELEK20121731,ZhaoHypergraphs}) offers even more complicated models. In those, the parameterizing function would not be $k$-dimensional, but rather $(2^k-2)$-dimensional, $U:\OurSpace^k\times \Lambda^{2^k-(k+2)}\to [0,\infty)$ (subject to certain symmetries). The right way to think about the power $2^k-(k+2)$ is that it represents all subsets of $[k]$ of sized $2,3,\ldots,k-1$. To create a random formula $\phi\sim\OriginalKSAT(n,U)$, we generate a random collection $\{x_S\in \Lambda\}_{S\subset [k],|S|=2,3,\ldots,k-1}$. The clause
$\{\mathfrak{s}_1 v_{i_1}, \ldots,\mathfrak{s}_k v_{i_k}\}$ is inserted into $\phi$ independently from other choices with probability
\[
\min\left\{1,
\frac{U\left((x_{\{i_1\}},\mathfrak{s}_1),\ldots,(x_{\{i_k\}},\mathfrak{s}_k),\big(x_S\big)_{S\subset [k],|S|=2,3,\ldots,k-1}\right)}{n^{k-1}}\;.
\right\}
\]
Given the difficulty of the homogeneous case of random $k$-SAT for $k\ge 3$ (\cite{MR4429261}), it is hopeless to obtain a complete picture of the inhomogeneous model. But there are some intermediate goals. For example, for a given $U$, it is plausible that there exists $c_U\in (0,+\infty]$ such that for $c<c_U$, $\OriginalKSAT(n,cU)$ is asymptotically almost surely satisfiable and for $c>c_U$, $\OriginalKSAT(n,cU)$ is asymptotically almost surely unsatisfiable.\footnote{The case $c_U=\infty$ corresponds to $U$s for which $\OriginalKSAT(n,cU)$ is almost asymptotically surely satisfiable for every $c>0$. This is a higher-dimensional counterpart to Proposition~\ref{prop:rhoNula}.}

\addcontentsline{toc}{section}{Index}
\printindex

\bibliographystyle{acm}
\addcontentsline{toc}{section}{Bibliography}
\bibliography{references.bib}

\section*{Contact details}
\noindent\begin{tabular}{ll} 
\emph{Postal address:} & Institute of Computer Science of the Czech Academy of Sciences \\ 
 & Pod Vodárenskou věží 2 \\ 
 & 182~00, Prague \\ 
 & Czechia\\
& \\
\emph{Email:}&\texttt{hladky@cs.cas.cz}, \texttt{savicky@cs.cas.cz}
\end{tabular}
\appendix

\end{document}